# A central-upwind scheme for open water flow in a wet/dry multiply-connected channel network

Sergii Kivva[a], Mark Zheleznyak[b], Alexander Pilipenko[a], Vasyl Yoschenko[b]

[a] Institute of Mathematical Machines & System Problems, Ukrainian National Academy of Sciences, 42, Prospect Glushkova, Kiev, 03680, Ukraine

[b] Institute of Environmental Radioactivity, Fukushima University, 1 Kanayagawa, Fukushima, 960-1296, Japan

**Abstract**. Our goal was to develop a robust algorithm for numerical simulation of one-dimensional shallow-water flow in a complex multiply-connected channel network with arbitrary geometry and variable topography. We apply a central-upwind scheme with a novel reconstruction of the open water surface in partially flooded cells that does not require additional correction. The proposed reconstruction and an exact integration of source terms for momentum conservation equation provide positivity preserving and well-balanced features of the scheme for various wet-dry states. We use two models based on continuity equation and mass and momentum conservation equations integrated for a control volume around the channel junction to it treatment. These junction models permit to simulate a subcritical and supercritical flow in a channel network. Numerous numerical experiments demonstrate the robustness of the proposed numerical algorithm and a good agreement of numerical results with exact solutions, experimental data, and results of the previous numerical studies. The proposed new specialized test on inundation and drying of initially dry channel network shows the merits of the new numerical algorithm to simulate the subcritical/supercritical open water flows in the network.

**Key words:** channel network with irregular geometry, shallow water equations, central-upwind scheme, well-balanced, positivity-preserving, wetting-drying

1. **Introduction**

The Saint Venant equations are common basis for the 1D modeling of the open flow in river and channel networks. Among different numerical methods in the various computer codes used in the hydraulic engineering (e.g. HEC-RAS [1], ISIS 1D [2], FEQ [3], CHARIMA [4], NETSTARS [5], FLDWAV [6] ) the Preissmann [7] scheme is one of the most widely used for the solution of one-dimensional unsteady open-flow channel problem. The Preissmann scheme is a weighted four-point nonlinear implicit finite-difference scheme. The advantage of this scheme is an efficient algorithm for the numerical solution of the Saint Venant equations in a multiply connected channel network. One of the disadvantages is that it does not provide positivity preserving of a numerical solution and therefore it cannot be applied to simulate water flow in drying channels. The positivity of water flow depth is especially important for schemes used for simulations of mountain rivers flow. Steep bottom slopes can generate artificial "drying" and following crash of the calculations due to the negative values of depth in numerical solutions if Preissmann scheme or other schemes without the positivity preserving are used for modeling.

Many difference schemes for hyperbolic equations that preserve non-negativity of water flow depth were developed after the pioneer publication of Godunov [8]. We refer the reader to [9], [10], [11], [12], [13], [14], [15], [16], [17], [18] where different approaches are presented. Among them, central [19] and central-upwind [20] schemes are widely used in a wide range of applications due to their simplicity, efficiency, and robustness. For the shallow water flow in open channels, these schemes were applied in [21], [22], [23]. Balbas and Karni [21] presented a central finite-volume scheme of second order accurate to simulate one-dimensional shallow water flows along channels with non-uniform rectangular cross-sections and bottom topography. In [22] Balbas and Hernandez-Duenas extended this scheme for channels with cross-sections of arbitrary shape and bottom topography. A central-upwind scheme with artificial viscosity was proposed by Hernandez-Duenas and Beljadid [23] for shallow water flows in channels with arbitrary geometry and variable topography.

In this paper, we implemented a central-upwind scheme [24] for open water flow in multiply-connected channel networks as a practical alternative to the Preissmann scheme. In the framework of the second order central-upwind schemes, we propose a novel reconstruction of water surface elevation that not requires additional correction in partially flooded cells. This reconstruction algorithm is a generalization of the reconstruction from [25].

The semi-discretization form of a central-upwind scheme is obtained by approximation of integral equations for the mass conservation and momentum balance. The integral form of momentum conservation equation includes terms that account for the forces due to changes in channel width and bed elevation. Assuming that the water and bed surface elevations and channel cross-section depend linearly on the water and bed elevations, and channel cross sections defined at cell faces, we can exactly calculate the source terms integrals on any internal cell interval. The exact values of the source terms integrals with our proposed reconstruction will guarantee us that the hydrostatic fluxes at a cell are balanced by the source terms that account for the sloping bed and variable channel width for various wetting/drying steady states at rest.

A treatment of channel junction in the modeling of open water flow in a multiply-connected channel network is one of the challenges. In many used hydraulic engineering codes, it is assumed that two groups of hydraulic conditions are held at a channel junction [26]. One is that the sum of the discharges has to be zero at the junction, and second is that the water levels or energy levels at the ends of linked channels (reaches) are equal at the junction. Such hydrodynamic models as CHARIMA, FEQ, MIKE-11 [27], RIVICE [28] and ONE-D [29] use the equality of water surface elevations at the junction. Others as HEC-RAS and FEQ are energy based models. Usually, these hydraulic conditions are specified as interior boundary conditions. In general case, the assumption of equality of inflow and outflow discharges at a junction does not permit in principle to simulate open water flow in a channel network under wetting/drying conditions. Yoshioka et al. [30, 31] use as the internal boundary conditions at the junction the continuity equation discretized on a dual mesh around it and the momentum flux at junction calculated as weighted linear combination of the momentum fluxes inflowing into it. Sanders et al. [32] and Bellamoli et al. [33] simulated open water flow in channel networks by nesting a 2D model at junctions. Contarino et al. [34] developed the implicit solver for the junction-generalized Riemann problem which was applied to construct a high-order ADER scheme for stiff hyperbolic balance laws in networks. We consider a channel junction as a computational cell that may have more than two inflows and outflows. For such cell, we integrate the continuity and momentum conservation equations over a control volume around the channel junction. We use in simulations two kinds of junction models: 1) based only on continuity equations for subcritical flow, and 2) based on mass and momentum conservation equations for supercritical flow. The first junction model is an analog of the equality of water surface elevations model.

Constraints on time step under which the scheme is positivity preserving may be very restrictive since they depend on the ratio of wetted cross-sectional areas at cell faces. To overcome this, similarly to [25], [35], we limit the outflow fluxes from draining cell by introducing the so-called local draining time that is smaller than the CFL time restriction. This approach provides positivity preserving of the scheme without a reduction of the CFL time restriction.

The paper is organized as follows. Channel network definition is presented in Section 2. Open flow water equations in channels are given in Section 3. In Section 4.1 we shortly describe the central-upwind scheme from [24] implemented for the Saint Venant equations. In Section 4.2, we present relations for exact discretization of source terms in the Saint Venant equations. New reconstruction method we describe in Section 4.3. Slope friction treatment and channel junction models are given in Section 4.4 and Section 4.5, respectively. Boundary conditions, which usually applied in modeling a subcritical open flow in a channel, are presented in Section 4.6. Positivity preserving and well-balancing properties of the scheme are proven in Sections 4.7 and 4.8, respectively. Section 5 contains different numerical tests illustrated the merits of the scheme. Section 6 provides some concluding remarks.

## 2. Channel network definition

The definition sketch of a channel network is shown in Fig.1 that is similar to CHARIMA [4]. Any channel network we will consider as a certain directed graph in which the edges are separate channels (*links*), and the vertices are points (*nodes*) in which channels may merge into a single channel or one channel may diverge into separate ones. The graph may also include loops that correspond to channels around islands.

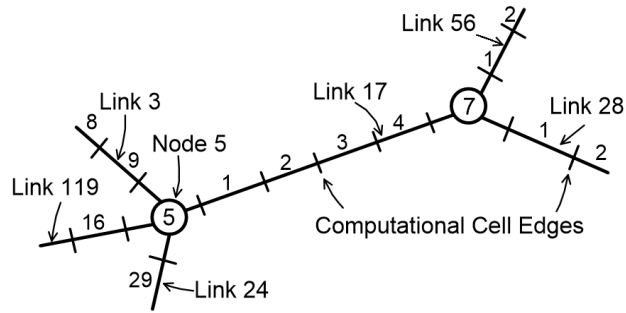

**Fig.1.** Definition sketch of the channel network.

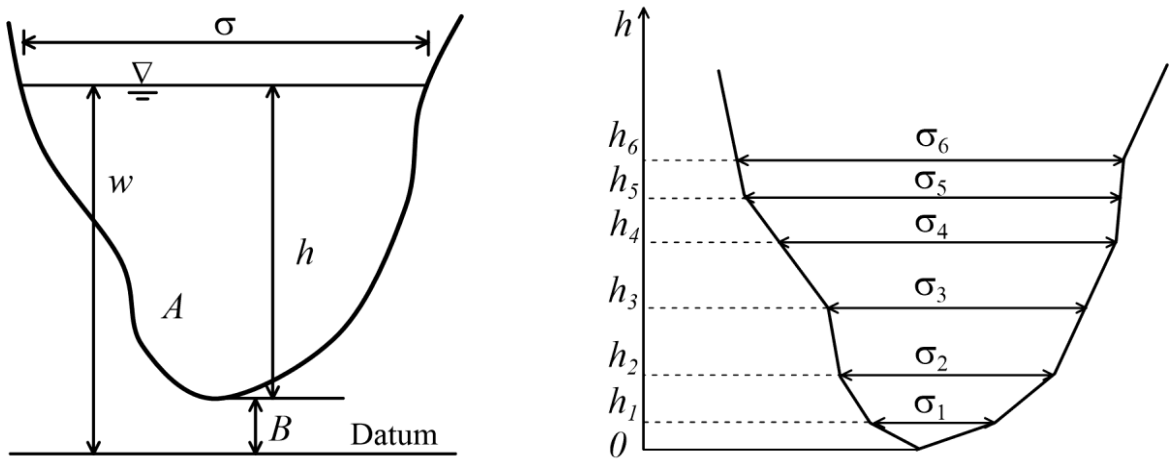

**Fig.2.** Channel cross-section and its piecewise linear approximation.

The *link* represents flow path between two *nodes* and each *link l* is divided into non-uniform computational cells, which are numbered from *1* at the upper end of the link to $n_l$ at the lower end of the link. Remark that "lower" and "upper" generally, although not necessarily, mean "downstream" and "upstream", respectively. Each link should have at least one computational cell.

Natural channel cross-section is replaced by its piecewise linear approximation (Fig.2) which assumed to be known for the each computational cell edge.

## 3. Channel water flow equations

The governing equations of one-dimensional shallow water flows for natural channel with irregular cross-section are given by the Saint Venant equations:

$$\frac{\partial A}{\partial t} + \frac{\partial Q}{\partial x} = 0 \qquad (1)$$

$$\frac{\partial Q}{\partial t} + \frac{\partial \left(Q^2/A + gI_1\right)}{\partial x} = gI_2 - gA(S_0 + S_f) \qquad (2)$$

where $t$ is the time; $x$ is the distance along the longitudinal axis of the channel; $Q$ is the discharge; $A$ is the wetted cross-section area; $g$ is the gravitational constant; $S_f$ is the friction slope due to the bed resistance; $S_0$ is the bed slope; $I_1$ is the hydrostatic pressure force; $I_2$ is the wall pressure force.

The friction slope $S_f$ is assumed to be given by the Manning equation

$$S_f = \frac{n^2 Q|Q|}{A^2 R^{4/3}} \qquad (3)$$

where $n$ is the Manning's roughness coefficient; $R$ is the cross-section hydraulic radius.

The hydrostatic pressure force, wall pressure force, wetted cross-section area and bed slope are defined, respectively, as

$$I_1 = \int_0^{h(x,t)} (h-y)\sigma(x,y)dy; \quad I_2 = \int_0^{h(x,t)} (h-y)\frac{\partial \sigma(x,y)}{\partial x}dy \qquad (4)$$

$$A = \int_0^{h(x,t)} \sigma(x,y)dy; \quad S_0 = \frac{\partial B}{\partial x} \qquad (5)$$

where $h$ is the water depth; $\sigma(x,h)$ is the channel width of a cross-section at point $x$ and water depth $h$; $B$ is the bed surface elevation (talweg or thalweg in geomorphology).

## 4. Numerical discretization

We consistently renumber all channels (links) and nodes in a network, that is, each channel and node have its own unique number. The numbering of channels does not depend on the numbering of nodes.

In general, index of the computational cells consists of the channel number and the cell number. The edges of the channel connected to one node are indexed by the node number and the channel number. Further, to simplify formulas where there will be no confusion, we omit some of the indices.

For simplicity, we assume that a channel geometry and a bottom topography do not change with time and depend only on the spatial coordinate. The other variables depend both on time and space. Also, where it will be possible, we will omit the temporal and spatial arguments.

*4.1. Semi-discrete central-upwind scheme*

A central-upwind scheme proposed by Kurganov and Petrova in [24] is used for numerical discretization of the Saint Venant equations (1)-(2).

We divide the spatial domain into the grid cells $\left[x_{j-1/2}, x_{j+1/2}\right]$ of the length $\Delta x_j$, where $x_j$ is the center of a grid cell, and denote by $\bar{U}_j(t)$ the cell averages of the solution $U = (A, Q)^T$ of (1)-(2) computed at time $t$

$$\bar{U}_j = \frac{1}{\Delta x_j} \int_{x_{j-1/2}}^{x_{j+1/2}} U(x,t) dx \tag{6}$$

Then, a semi-discretization of (1)-(2) can be written as the following system of ODEs

$$\frac{d}{dt} \bar{U}_j(t) = -\frac{H_{j+1/2}(t) - H_{j-1/2}(t)}{\Delta x_j} + \frac{1}{\Delta x_j} \int_{x_{j-1/2}}^{x_{j+1/2}} [S(U,x) + G(U,x)] dx \tag{7}$$

where $H_{j+1/2}$ are numerical fluxes at the cell interfaces $x_{j\pm 1/2}$, $S(U,x) = (0, gI_2 - gAS_0)^T$ and $G(U,x) = (0, -gAS_f)^T$.

The central-upwind numerical fluxes $H_{j+1/2}$ are given by

$$H_{j+1/2} = \frac{a_{j+1/2}^+ F(U_{j+1/2}^-, B_{j+1/2}) - a_{j+1/2}^- F(U_{j+1/2}^+, B_{j+1/2})}{a_{j+1/2}^+ - a_{j+1/2}^-} + \frac{a_{j+1/2}^+ a_{j+1/2}^-}{a_{j+1/2}^+ - a_{j+1/2}^-} \left[ U_{j+1/2}^+ - U_{j+1/2}^- \right] \tag{8}$$

where $F(U,B) = (Q, Q^2/A + gI_1)^T$ is the flux at the cell interface, and the one-sided local speeds are obtained using the eigenvalues of the Jacobian

$$a_{j+1/2}^+ = \max\left\{0, u_{j+1/2}^+ + c_{j+1/2}^+, u_{j+1/2}^- + c_{j+1/2}^-\right\} \tag{9}$$

$$a_{j+1/2}^- = \min\left\{0, u_{j+1/2}^+ - c_{j+1/2}^+, u_{j+1/2}^- - c_{j+1/2}^-\right\}$$

Here $u = Q/A$, $c = \sqrt{gA/\sigma_T}$ and $\sigma_T$ is the width of the channel at the water surface. The wetted water area at the cell interface may be very small and may lead to large values of flow velocity. In order to prevent this, we use the regularization technique suggested in [24]

$$u = \frac{\sqrt{2}AQ}{\sqrt{A^4 + \max(A^4, \varepsilon)}} \tag{10}$$

where $\varepsilon$ is a small apriori chosen positive number. For consistency, the values of the discharge at cell interfaces are recomputed using $Q^{\pm}_{j\pm1/2} = A^{\pm}_{j\pm1/2} u^{\pm}_{j\pm1/2}$, where $u^{\pm}_{j\pm1/2}$ are given by (10).

The interface values $U^{\pm}_{j+1/2}$ are obtained from the cell averages $\bar{U}_j$ by a piecewise linear reconstruction which will be described below.

### 4.2. Auxiliary relationships

In this section, we obtain relationships for computing the water volume between two specified channel cross-sections for an arbitrary free surface of the water. We also obtain formulas for calculating the hydrostatic pressure force $I_1$ and the wall pressure force $I_2$.

We replace the bed function $B$ and wetted area $A$ on the $\left[x_{j-1/2}, x_{j+1/2}\right]$ by their piecewise linear approximation

$$B(x) = B_{j-1/2} + (B_{j+1/2} - B_{j-1/2}) \frac{x - x_j}{\Delta x_j}, \tag{11}$$

$$A(x,t) = \frac{1}{\Delta x_j} \int_0^{h(x,t)} \left[\sigma_{j+1/2}(y)(x - x_{i-1/2}) + \sigma_{i-1/2}(y)(x_{j+1/2} - x)\right] dy \tag{12}$$

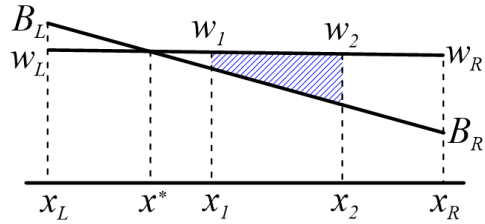

**Fig.3.** Water surface $w$ in a partially flooded cell $[x_L, x_R]$. Calculation of the water volume between $x_1$ and $x_2$.

Assuming that the open water surface $w$ is linear in a computational cell $[x_L, x_R]$, below we will use the following relations. Water volume $V(x_1, x_2; w_L, w_R)$ between points $x_1$ and $x_2$ (Fig.3.) equals to

$$V(x_1, x_2; w_L, w_R) = \int_{x_1}^{x_2} A(x,h) dx = \int_{x_1}^{x_2} \int_0^{h(x,t)} \sigma(x,y) dy dx$$

Changing the order of integration, we obtain

$V(x_1, x_2; w_L, w_R) =$

$$= \frac{(x_2 - x_1)(2x_R - x_1 - x_2)}{2\Delta x} \int_0^{h_1} \sigma_L dy - \frac{(x_R - x_2)^2}{2\Delta x} \int_{h_1}^{h_2} \sigma_L dy + \frac{\Delta x}{2\Delta h^2} \int_{h_1}^{h_2} \sigma_L (h_R - y)^2 dy + \tag{13}$$

$$+ \frac{(x_2 - x_1)(x_1 + x_2 - 2x_L)}{2\Delta x} \int_0^{h_1} \sigma_R dy + \frac{(x_2 - x_L)^2}{2\Delta x} \int_{h_1}^{h_2} \sigma_R dy - \frac{\Delta x}{2\Delta h^2} \int_{h_1}^{h_2} \sigma_R (y - h_L)^2 dy$$

where $h_i = w_i - B_i$, $\Delta x = x_R - x_L$ and $\Delta h = h_R - h_L$. Note that for $h_1 = h_2$, the second, third, fifth and sixth terms on the right-hand side of (13) are zero.

Similarly, the hydrostatic pressure force $I_1$ at a point $x_1$ can be calculated as

$$I_1(x_1, h_1) = \frac{(x_R - x_1)}{\Delta x} \int_0^{h_1} (h_1 - y)\sigma_L(y)dy + \frac{(x_1 - x_L)}{\Delta x} \int_0^{h_1} (h_1 - y)\sigma_R(y)dy \tag{14}$$

Integration of the wall pressure force $I_2$ and $AS_0$ over an interval $[x_1, x_2]$ yields the following expressions

$$I_2(x_1, x_2) = \int_{x_1}^{x_2} I_2 dx = \frac{x_2 - x_1}{\Delta x} \int_0^{h_1} (\sigma_R - \sigma_L)\left(\frac{h_1 + h_2}{2} - y\right)dy + \frac{1}{2\Delta h} \int_{h_1}^{h_2} (\sigma_R - \sigma_L)(h_2 - y)^2 dy \tag{15}$$

$$B_x(x_1, x_2) = \int_{x_1}^{x_2} A\frac{\partial B}{\partial x} dx = \frac{B_R - B_L}{2\Delta x} V(x_1, x_2; w_L, w_R) \tag{16}$$

*Remark.*

1. Note that for linear approximation of the water surface and bottom on an interval $\left[x_{j-1/2}, x_{j+1/2}\right]$, if the values of $w_{j\mp1/2}^{\pm}$ and $B_{j\pm1/2}$ are known at the cell faces $x_{j\mp1/2}$, then it is easy to calculate $x^*$ and integrals (13)-(16) for any internal interval $[x_1, x_2]$ and nonnegative $h(x,t)$.

2. The integrals (13)-(16) can be also computed for any piecewise-linear approximation of the water surface on an interval $\left[x_{j-1/2}, x_{j+1/2}\right]$.

### 4.3. Reconstruction of the interface values $U_{j+1/2}^{\pm}$

We denote by $h_{av,j}$, $w_j$ and $w_{st,j}$ the depth of water which surface is parallel to the cell bed, water surface elevation and the still water surface elevation in the cell $j$, respectively. In addition, we denote $h^{\geq 0} = \max(0, h)$.

Substituting $w_1 = w_2 = w_{st,j}$ into equation (13) and taking into account that $h_{j+1/2} - h_{j-1/2} = B_{j-1/2} - B_{j+1/2} = -\Delta B$, we obtain the following equation for finding $w_{st,j}$

$$\frac{1}{2}\int_0^{(w_{st,j} - B_{j-1/2})^{\geq 0}} \sigma_{j-1/2} dy + \frac{1}{2\Delta B^2}\int_{(w_{st,j} - B_{j-1/2})^{\geq 0}}^{(w_{st,j} - B_{j+1/2})^{\geq 0}} (w_{st,j} - B_{j+1/2} - y)^2 \sigma_{j-1/2} dy +$$

$$+ \frac{1}{2}\int_0^{(w_{st,j} - B_{j+1/2})^{\geq 0}} \sigma_{j+1/2} dy - \frac{1}{2\Delta B^2}\int_{(w_{st,j} - B_{j-1/2})^{\geq 0}}^{(w_{st,j} - B_{j+1/2})^{\geq 0}} (y - w_{st,j} + B_{j-1/2})^2 \sigma_{j+1/2} dy = \overline{A}_j \tag{17}$$

Equation (17) is a quartic equation for the piecewise linear approximation of channel cross-section and we can use the Newton-Raphson method or Ferrari's method for its solution. We calculate $w_{st,j}$ only for cells in which still water may occur.

Similarly, substituting $h_1 = h_2 = h_{av,j}$ into equation (13), we obtain the following equation for finding $h_{av,j}$

$$\frac{1}{2}\int_0^{h_{av,j}} \sigma_{j-1/2} dy + \frac{1}{2}\int_0^{h_{av,j}} \sigma_{j+1/2} dy = \bar{A}_j \qquad (18)$$

which is a quadratic equation since

$$\int_0^h \sigma_{j+1/2} dy = \sum_{k=0}^{p-1}(h_{j+1/2,k+1} - h_{j+1/2,k})\frac{\sigma_{j+1/2,k} + \sigma_{j+1/2,k+1}}{2} + \qquad (19)$$

$$+(h - h_{j+1/2,p})\left[\sigma_{j+1/2,p} + \frac{1}{2}\frac{\sigma_{j+1/2,p+1} - \sigma_{j+1/2,p}}{h_{j+1/2,p+1} - h_{j+1/2,p}}(h - h_{j+1/2,p})\right], \qquad \text{for } h_{j+1/2,p} < h \le h_{j+1/2,p+1}$$

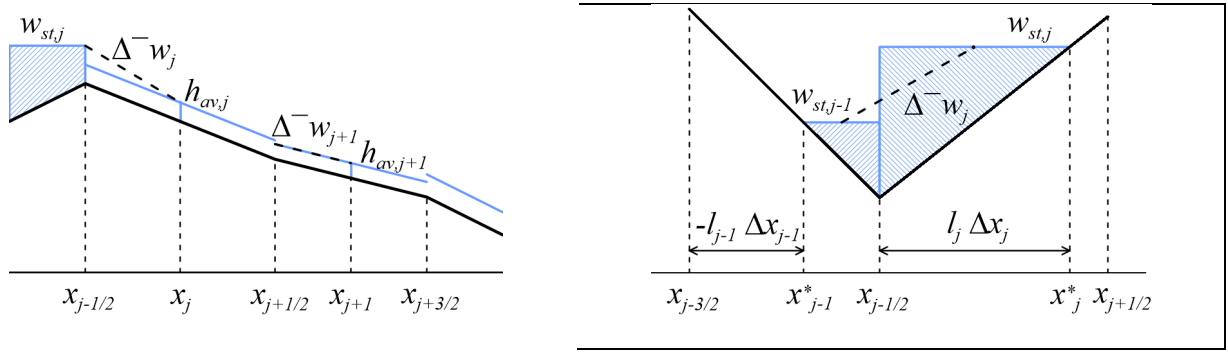

**Fig.4.** Approximation of the backward difference derivative $\Delta^- w$.

We divide conditionally all computational cells into two groups: "wet" and "dry". A cell $j$ will be called "wet" if $w_{st,j} \ge \max(B_{j-1/2}, B_{j+1/2})$, and "dry" otherwise. We assume $w_j = w_{st,j}$ for the "wet" cells and $w_j = h_{av,j} + (B_{j+1/2} + B_{j-1/2})/2$ for the "dry" cell. For "dry" cells in which still water may occur, we assume $w_j = w_{st,j}$ and $l_j = (w_{st,j} - B_{j-1/2})/|\Delta B_j|$. For all other cells $l_j = 1$.

Depending on combination of two neighboring cells, we use different formulas for approximating the forward and backward difference derivatives for $w$ and $Q$. We approximate backward difference derivative $\Delta^-$ by the following way (Fig.4):

1. $j-1$ and $j$ cell are both "wet". Then

$$\Delta^- w_j = \frac{w_{st,j} - w_{st,j-1}}{x_j - x_{j-1}}; \quad \Delta^- Q_j = \frac{\bar{Q}_j - \bar{Q}_{j-1}}{x_j - x_{j-1}} \qquad (20)$$

2. $j-1$ cell is "wet" and $j$ cell is "dry"

- if $B_{j+1/2} > B_{j-1/2}$

$$\Delta^- w_j = 2\frac{w_{st,j} - w_{st,j-1}}{\Delta x_{j-1} + l_j \Delta x_j}; \qquad \Delta^- Q_j = 2\frac{\bar{Q}_j - \bar{Q}_{j-1}}{\Delta x_{j-1} + l_j \Delta x_j} \qquad (21)$$

- otherwise

$$\Delta^- w_j = 2\frac{B_{j-1/2} + h_{av,j} - w_{st,j-1}}{\Delta x_j}; \qquad \Delta^- Q_j = 2\frac{\bar{Q}_j - \bar{Q}_{j-1}}{\Delta x_j} \tag{22}$$

3. *j-1* cell is "dry" and *j* cell is "wet"

- if $B_{j-3/2} > B_{j-1/2}$

$$\Delta^- w_j = 2\frac{w_{st,j} - w_{st,j-1}}{(1+l_{j-1})\Delta x_{j-1} + \Delta x_j}; \qquad \Delta^- Q_j = 2\frac{\bar{Q}_j - \bar{Q}_{j-1}}{(1+l_{j-1})\Delta x_{j-1} + \Delta x_j} \tag{23}$$

- otherwise

$$\Delta^- w_j = 2\frac{w_{st,j} - B_{j-1/2} - h_{av,j-1}}{\Delta x_j}; \qquad \Delta^- Q_j = 2\frac{\bar{Q}_j - \bar{Q}_{j-1}}{\Delta x_j} \tag{24}$$

4. *j-1* and *j* cells are both "dry"

- if $B_{j-3/2} > B_{j-1/2}$ and $B_{j+1/2} > B_{j-1/2}$

$$\Delta^- w_j = 2\frac{w_{st,j} - w_{st,j-1}}{(1+l_{j-1})\Delta x_{j-1} + l_j \Delta x_j}; \qquad \Delta^- Q_j = 2\frac{\bar{Q}_j - \bar{Q}_{j-1}}{(1+l_{j-1})\Delta x_{j-1} + l_j \Delta x_j} \tag{25}$$

- otherwise

$$\Delta^- w_j = \frac{B_{j+1/2} - B_{j-1/2}}{\Delta x_j}; \qquad \Delta^- Q_j = 0 \tag{26}$$

The forward derivatives $\Delta^+$ for *w* and *Q* are calculated in the same way.

The interface point-values $w^{\pm}_{j\mp 1/2}$ are obtained by a piecewise linear reconstruction

$$w^+_{j-1/2} = \begin{cases} w_j - w_{x,j}\dfrac{l_j}{2}\Delta x_j, & \text{if } l_j \geq 0 \\ w_j - w_{x,j}\dfrac{1-l_j}{2}\Delta x_j, & \text{otherwise} \end{cases} \tag{27}$$

and

$$w^-_{j+1/2} = \begin{cases} w_j + w_{x,j}(1-\dfrac{l_j}{2})\Delta x_j, & \text{if } l_j \geq 0 \\ w_j + w_{x,j}\dfrac{1+l_j}{2}\Delta x_j, & \text{otherwise} \end{cases} \tag{28}$$

According to [19] this reconstruction will be second-order accurate if the approximate values of the derivatives $w_{x,j}$ are at least first-order approximations of corresponding exact derivatives. To ensure the non-oscillatory property in our numerical scheme we evaluate $w_{x,j}$ using the minmod limiter [13], [36], [37]

$$w_{x,j} = \min\bmod(\Delta^- w_j, \Delta^+ w_j) \tag{29}$$

where $\min\bmod(a,b) = \dfrac{\mathrm{sgn}(a)+\mathrm{sgn}(b)}{2}\min(|a|,|b|)$.

The interface values $Q^{\pm}_{j\mp 1/2}$ we calculate in the same way. In the next step, water depth $h^{\pm}_{j\pm 1/2}$ at the cell faces is computed from the water surface elevation. Then, the area $A^{\pm}_{j\pm 1/2}$ is calculated from water depth based on channel geometry. Finally, we determine the numerical fluxes $H_{j\pm 1/2}$ and $\int_{x_{j-1/2}}^{x_{j+1/2}} S(U,x)dx$ from (7), (8) using formulas (14)-(16) to compute $I_1$, $I_2(x_{j-1/2}, x_{j+1/2})$, and $B_x(x_{j-1/2}, x_{j+1/2})$.

### 4.4. Friction slope treatment and time evolution

We calculate the second component of $\int_{x_{j-1/2}}^{x_{j+1/2}} G(U,x)dx$ using the midpoint rule and de-singularization procedure

$$\frac{1}{\Delta x_j}\int_{x_{j-1/2}}^{x_{j+1/2}} G(U,x)dx = gn^2 \frac{|\bar{Q}_j|}{\max(\bar{A}\,R^{4/3},\varepsilon)} \bar{Q}_j = G^{(2)}_j \bar{Q}_j \tag{30}$$

In the general case, the friction slope term in (7) can be a source of stiffness for the ODE system when the depth of water is very small [38]. The explicit treatment of the friction term imposes a severe time step restriction, as the result of which time step should be in several times less than a typical time step under the CFL conditions. To overcome this difficulty, we use the forward Euler method of time integration of the ODE system (7) with an implicit treatment of only a part of the friction slope term. In result, we obtained the following fully discrete central-upwind scheme

$$\bar{A}_j(t+\Delta t) = \bar{A}_j(t) - \frac{\Delta t}{\Delta x_j}\left(H^{(1)}_{j+1/2}(t) - H^{(1)}_{j-1/2}(t)\right) \tag{31}$$

$$\bar{Q}_j(t+\Delta t) = \frac{\bar{Q}_j(t) - \dfrac{\Delta t}{\Delta x_j}\left(H^{(2)}_{j+1/2}(t) - H^{(2)}_{j-1/2}(t) - g I_2(x_{j-1/2}, x_{j+1/2}) + gB_x(x_{j-1/2}, x_{j+1/2})\right)}{1+\Delta t\, G^{(2)}_j(t)} \tag{32}$$

*Remark.* Note that the scheme (31)-(32) is first-order accurate in time. A higher-order time discretization can be obtained using semi-implicit solvers proposed by Chertok et al. [39] which are based on the modification of explicit SSP Runge-Kutta methods [40]. The second-order semi-implicit method for a system ODE with a non-stiff term $f(u,t)$ and stiff damping term $g(u,t)u$ such that $g(u,t)$ is a diagonal matrix with non-positive elements

$$u' = f(u,t) + g(u,t)u$$

can be written as

$$u^{(1)} = (E - \Delta t\, g^n)^{-1}(u^n + \Delta t\, f^n)$$

$$u^{(2)} = \frac{1}{2}u^n + \frac{1}{2}(E - \Delta t\, g^{(1)})^{-1}(u^{(1)} + \Delta t\, f^{(1)})$$

$$u^{n+1} = \left[E + (\Delta t g^{(2)})^2\right]^{-1}\left[u^{(2)} - (\Delta t)^2 f^{(2)} g^{(2)}\right]$$

### *4.5. Node of junction*

For modeling of the hydraulic conditions at a junction, we will use two approaches: (1) assuming that only mass conservation equation is hold at the junction, and (2) based on mass and momentum conservations.

Consider a node $s$ that connects two and more channels, $J_{s,in}$ is a set of channels which inflow and $J_{s,out}$ that outflow from it (Fig.5). Denote $x^{in}_{s,i}$ and $x^{out}_{s,j}$ coordinates of the ends of channels connected to the node $s$. Channel cross-sections and bottom levels are known at these points. In the first case, water surface elevation $Y_s$, or water surface elevation $Y_s$ and water discharge $Q_s$ for the second case at the junction are also given.

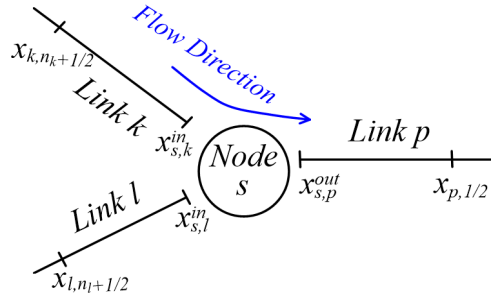

**Fig.5**. Sketch of a node that is the junction of channels.

Assuming that water surface is horizontal, the continuity equation over the control volume around a junction $s$ can be written as follow

$$F_s^A\left(Y_s(t+\Delta t)\right) = \sum_{i \in J_{s,in}}\left[V(x_{i,n_i+1/2}, x^{in}_{s,i}; Y_s(t+\Delta t), Y_s(t+\Delta t)) - V(x_{i,n_i+1/2}, x^{in}_{s,i}; Y_s(t), Y_s(t))\right] +$$

$$+ \sum_{j \in J_{s,out}}\left[V(x^{out}_{s,j}, x_{j,1/2}; Y_s(t+\Delta t), Y_s(t+\Delta t)) - V(x^{out}_{s,j}, x_{j,1/2}; Y_s(t), Y_s(t))\right] - \qquad (33)$$

$$-\Delta t \sum_{i \in J_{s,in}}\left\{\frac{a^+_{i,n_i+1/2}Q^-_{i,n_i+1/2} - a^-_{i,n_i+1/2}Q^+_{i,n_i+1/2}}{a^+_{i,n_i+1/2} - a^-_{i,n_i+1/2}} + \frac{a^+_{i,n_i+1/2}a^-_{i,n_i+1/2}}{a^+_{i,n_i+1/2} - a^-_{i,n_i+1/2}}\left[A^+_{i,n_i+1/2} - A^-_{i,n_i+1/2}\right]\right\} +$$

$$+\Delta t \sum_{j \in J_{s,out}}\left\{\frac{a^+_{j,1/2}Q^-_{j,1/2} - a^-_{j,1/2}Q^+_{j,1/2}}{a^+_{j,1/2} - a^-_{j,1/2}} + \frac{a^+_{j,1/2}a^-_{j,1/2}}{a^+_{j,1/2} - a^-_{j,1/2}}\left[A^+_{j,1/2} - A^-_{j,1/2}\right]\right\} = 0$$

Accordingly, the equation of momentum conservation can be written in the form

$$F_s^Q\left(\bar{Q}_s(t+\Delta t)\right) = \left(\bar{Q}_s(t+\Delta t) - \bar{Q}_s(t)\right)\left[\sum_{i\in J_{s,in}}(x_{s,i}^{in} - x_{i,n_i+1/2}) + \sum_{j\in J_{s,out}}(x_{j,1/2} - x_{s,j}^{out})\right] +$$

$$-\Delta t \sum_{i\in J_{s,in}}\left\{H_{i,n_i+1/2}^{(2)} - gI_1(x_{s,i}^{in}, h_s) + gI_2(x_{x_{i,n_i+1/2}}, x_{s,i}^{in}) - gB_x(x_{x_{i,n_i+1/2}}, x_{s,i}^{in})\right\} + \quad (34)$$

$$+\Delta t \sum_{j\in J_{s,out}}\left\{H_{j,1/2}^{(2)} - gI_1(x_{s,j}^{out}, h_s) - gI_2(x_{s,j}^{out}, x_{j,1/2}) + gB_x(x_{s,j}^{out}, x_{j,1/2})\right\} -$$

$$-\Delta t \bar{Q}_s(t+\Delta t)\left[\sum_{i\in J_{s,in}}G_{s,i}^{(2)}(x_{s,i}^{in} - x_{i,n_i+1/2}) + \sum_{j\in J_{s,out}}G_{s,j}^{(2)}(x_{j,1/2} - x_{s,j}^{out})\right] = 0$$

The interface values $A_{j+1/2}^{\pm}$ are calculated from water surface elevation, which is reconstructed in a similar way as described in section 4.3. The interface values $Q_{j+1/2}^{\pm}$ are computed: in the first case as $Q_{i,n_i+1/2}^{+} = Q_{i,n_i+1/2}^{-}$, $Q_{j,1/2}^{-} = Q_{j,1/2}^{+}$, and in the second case as $Q_{i,n_i+1/2}^{+} = Q_{j,1/2}^{-} = \bar{Q}_s$.

The water surface elevation $Y_s(t+\Delta t)$ and the water discharge $Q_s(t+\Delta t)$ are found from the solution of equations (33)-(34). The Newton-Raphson method is used to solve the nonlinear equation (33). Note that the derivative of $F_s^A(Y_s(t+\Delta t))$ with respect to $Y_s(t+\Delta t)$ can be expressed as

$$\frac{\partial F_s^A}{\partial Y_s} = \frac{1}{2}\sum_{i\in J_{s,in}}\left(x_{s,i}^{in} - x_{i,n_i+1/2}\right)\left[\sigma_{i,n_i+1/2}\left(h_{i,n_i+1/2}^{\geq 0}\right) + \sigma_{s,i}\left(h_{s,i}^{\geq 0}\right) - \right.$$

$$-\left(B_{s,i} - B_{i,n_i+1/2}\right)^{-2}\left\{\left[\sigma_{i,n_i+1/2}\left(h_{i,n_i+1/2}^{\geq 0}\right) + \sigma_{s,i}\left(h_{s,i}^{\geq 0}\right)\right]\left(B_{s,i} - B_{i,n_i+1/2}\right)^2 - \right.$$

$$\left.\left.-2\int_{(Y_s-B_{i,n_i+1/2})^{\geq 0}}^{(Y_s-B_{s,i})^{\geq 0}}\left[(Y_s - B_{s,i} - y)\sigma_{i,n_i+1/2}(y) + (y - Y_s + B_{i,n_i+1/2})\sigma_{s,i}(y)\right]dy\right\}\right] - \quad (35)$$

$$+\frac{1}{2}\sum_{j\in J_{s,out}}\left(x_{j,1/2} - x_{s,j}^{out}\right)\left[\sigma_{s,j}\left(h_{s,j}^{\geq 0}\right) + \sigma_{j,1/2}\left(h_{j,1/2}^{\geq 0}\right) - \right.$$

$$-\left(B_{j,1/2} - B_{s,j}\right)^{-2}\left\{\left[\sigma_{s,j}\left(h_{s,j}^{\geq 0}\right) + \sigma_{j,1/2}\left(h_{j,1/2}^{\geq 0}\right)\right]\left(B_{j,1/2} - B_{s,j}\right)^2 - \right.$$

$$\left.\left.-2\int_{(Y_s-B_{s,j})^{\geq 0}}^{(Y_s-B_{j,1/2})^{\geq 0}}\left[(Y_s - B_{j,1/2} - y)\sigma_{s,j}(y) + (y - Y_s + B_{s,j})\sigma_{j,1/2}(y)\right]dy\right\}\right]$$

where $h^{\geq 0} = \max(0,h)$. Note that for $h_{i,n_i+1/2} = h_{s,i}$ and $h_{j,1/2} = h_{s,j}$ the expressions in the corresponding curly brackets on the right-hand side of (35) are zero.

### *4.6. Boundary conditions*

For the numerical solution of the ODEs (7), we should specify values of the central-upwind numerical fluxes $H_{j+1/2}$ at an upstream and if required downstream boundaries of the computational area. We consider the most commonly used boundary conditions in the modeling of a subcritical fluid flow with a free surface in channel networks. In the following, we discuss the boundary treatment at the left boundary. Boundary treatment at the right boundary can be done by similar way.

*Discharge boundary conditions.* Let a water discharge $Q(t)$ has been provided on the left boundary. We will consider the following approximation of the boundary numerical fluxes $H_{1/2}$

$$Q = \frac{a_{1/2}^+ Q_{1/2}^- - a_{1/2}^- Q_{1/2}^+}{a_{1/2}^+ - a_{1/2}^-} + \frac{a_{1/2}^- a_{1/2}^+}{a_{1/2}^+ - a_{1/2}^-}\left(A_{1/2}^+ - A_{1/2}^-\right) \tag{36}$$

$$\frac{Q^2}{A_{1/2}^-} + gI_1(h_{1/2}^-) = \frac{a_{1/2}^+ \dfrac{(Q_{1/2}^-)^2}{A_{1/2}^-} - a_{1/2}^- \dfrac{(Q_{1/2}^+)^2}{A_{1/2}^+}}{a_{1/2}^+ - a_{1/2}^-} + \frac{a_{1/2}^- a_{1/2}^+}{a_{1/2}^+ - a_{1/2}^-}\left(Q_{1/2}^+ - Q_{1/2}^-\right) + g\frac{a_{1/2}^+ I_1(h_{1/2}^-) - a_{1/2}^- I_1(h_{1/2}^+)}{a_{1/2}^+ - a_{1/2}^-}$$

The values of $Q_{1/2}^-$ и $h_{1/2}^-$ are found from a solution of the system of equations (36).

*Surface boundary condition.* Let a water surface elevation $Y(t)$ is specified at the left boundary. Then we will use the following approximation of the boundary fluxes

$$Q_{1/2}^- = \frac{a_{1/2}^+ Q_{1/2}^- - a_{1/2}^- Q_{1/2}^+}{a_{1/2}^+ - a_{1/2}^-} + \frac{a_{1/2}^- a_{1/2}^+}{a_{1/2}^+ - a_{1/2}^-}\left(A_{1/2}^+ - A_{1/2}^-\right) \tag{37}$$

$$\frac{(Q_{1/2}^-)^2}{A(h)} + gI_1(h) = \frac{a_{1/2}^+ \dfrac{(Q_{1/2}^-)^2}{A_{1/2}^-} - a_{1/2}^- \dfrac{(Q_{1/2}^+)^2}{A_{1/2}^+}}{a_{1/2}^+ - a_{1/2}^-} + \frac{a_{1/2}^- a_{1/2}^+}{a_{1/2}^+ - a_{1/2}^-}\left(Q_{1/2}^+ - Q_{1/2}^-\right) + g\frac{a_{1/2}^+ I_1(h_{1/2}^-) - a_{1/2}^- I_1(h_{1/2}^+)}{a_{1/2}^+ - a_{1/2}^-}$$

where $h = Y - B_{1/2}$; $Q_{1/2}^-$ и $h_{1/2}^-$ are found from a solution of the equations (37).

We will also consider a more simple approximation of the surface boundary conditions. In this case, we assume that

$$Q_{1/2}^- = Q_{1/2}^+ \tag{38}$$

$$h_{1/2}^- = h$$

*Outflow boundary condition.* For outflow, we use the following extrapolations on the boundary

$$Q_{1/2}^- = Q_{1/2}^+ \tag{39}$$

$$h_{1/2}^- = h_{1/2}^+$$

### 4.7. Positivity preserving

In this section, we show that the resulting central-upwind scheme is positivity preserving. A sufficient condition for this is to ensure that at each time step no more water outflow from a cell than it is at the moment. For the positivity, our result is similar to the one obtained in [22].

For any computational cell $k$ of the link $p$, we define $\Delta t_{p,k}$ as

$$\Delta t_{p,k} = \min\left\{\frac{\Delta x_{p,k}}{a_{p,k-1/2}^+ - a_{p,k+1/2}^-}, \frac{\Delta x_{p,k}\overline{A}_{p,k}^n}{a_{p,k+1/2}^+ A_{p,k+1/2}^- - a_{p,k-1/2}^- A_{p,k-1/2}^+}\right\} \tag{40}$$

For any node s, which has two or more links, we denote by $\Delta t_s$

$$\Delta t_s = \min\left\{\min_i\left[\frac{\Delta x_{s,i}^{in}}{a_{i,n_i+1/2}^+}, \frac{-\Delta x_{s,i}^{in}\overline{A}_{s,i}^{in,n}}{a_{i,n_i+1/2}^- A_{i,n_i+1/2}^+}\right], \min_j\left[\frac{-\Delta x_{s,j}^{out}}{a_{j,1/2}^-}, \frac{\Delta x_{s,j}^{out}\overline{A}_{s,j}^{out,n}}{a_{j,1/2}^+ A_{j,1/2}^-}\right]\right\} \tag{41}$$

where $V(x_{s,j}^{out}, x_{j,1/2}; Y_s(t), Y_s(t)) = \Delta x_{s,j}^{out}\overline{A}_{s,j}^{out,n}$ and $\Delta x_{s,j}^{out} = x_{j,1/2} - x_{s,j}^{out}$. The values of $\overline{A}_{s,i}^{in,n}$ and $\Delta x_{s,i}^{in}$ are defined in a similar way.

**Theorem 1.** *Consider the semi-discrete central-upwind scheme (7)-(10), (33) with the piecewise linear reconstruction described in section 4.3 and the discretization of the source terms (14)-(16). Assume that the system of ODEs (7) is solved by the forward Euler method with an implicit treatment of the friction slope and for all computational cells and nodes $h_j^n \geq 0$. Then all $h_j^{n+1} \geq 0$, if*

$$\Delta t \leq \min\left\{\min_{p,k}\Delta t_{p,k}, \Delta t_s\right\} \tag{42}$$

*where $\Delta t_{p,k}$ and $\Delta t_s$ are calculated from (40)-(41).*

*Proof.* The cross-sectional area $A(h)$ is a nonnegative increasing function of water depth. Therefore, our task is to obtain conditions that will ensure the non-negativity of the cross-section area at the next time step.

Note, as follows from (9) that $a_{j+1/2}^+ \geq 0$, $a_{j+1/2}^- \leq 0$, $a_{j+1/2}^+ - u_{j+1/2}^+ \geq 0$, and $u_{j+1/2}^- - a_{j+1/2}^- \geq 0$ for all $j$. Moreover, the following inequalities $0 \leq \frac{u_{j+1/2}^- - a_{j+1/2}^-}{a_{j+1/2}^+ - a_{j+1/2}^-} \leq 1$, $0 \leq \frac{a_{j-1/2}^+ - u_{j-1/2}^+}{a_{j-1/2}^+ - a_{j-1/2}^-} \leq 1$ are satisfied for all $j$. Then for any computational cell $k$ of link $p$ from (31), we have

$$\overline{A}_{p,k}^{n+1} = \overline{A}_{p,k}^n - \frac{\Delta t}{\Delta x_{p,k}}a_{p,k+1/2}^+\left(\frac{u_{p,k+1/2}^- - a_{p,k+1/2}^-}{a_{p,k+1/2}^+ - a_{p,k+1/2}^-}\right)A_{p,k+1/2}^- - \frac{\Delta t}{\Delta x_{p,k}}a_{p,k+1/2}^-\left(\frac{a_{p,k+1/2}^+ - u_{p,k+1/2}^+}{a_{p,k+1/2}^+ - a_{p,k+1/2}^-}\right)A_{p,k+1/2}^+ +$$

$$+\frac{\Delta t}{\Delta x_{p,k}}a_{p,k-1/2}^+\left(\frac{u_{p,k-1/2}^- - a_{p,k-1/2}^-}{a_{p,k-1/2}^+ - a_{p,k-1/2}^-}\right)A_{p,k-1/2}^- + \frac{\Delta t}{\Delta x_{p,k}}a_{p,k-1/2}^-\left(\frac{a_{p,k-1/2}^+ - u_{p,k-1/2}^+}{a_{p,k-1/2}^+ - a_{p,k-1/2}^-}\right)A_{p,k-1/2}^+ \tag{43}$$

The third and fourth terms on the right-hand side of (43) are nonnegative so that $\overline{A}_{p,k}^{n+1} \geq 0$ it is sufficient if the following inequality is satisfied

$$\Delta x_{p,k}\overline{A}_{p,k}^n - \Delta t\left(a_{p,k+1/2}^+ A_{p,k+1/2}^- - a_{p,k-1/2}^- A_{p,k-1/2}^+\right) \geq 0 \tag{44}$$

For any node $s$, we rewrite (33) in the following way

$$\sum_{i\in J_{s,in}}\Delta x_{s,i}^{in}\overline{A}_{s,i}^{in,n+1} + \sum_{j\in J_{s,out}}\Delta x_{s,j}^{out}\overline{A}_{s,j}^{out,n+1} = \sum_{i\in J_{s,in}}\Delta x_{s,i}^{in}\overline{A}_{s,i}^{in,n} + \sum_{j\in J_{s,out}}\Delta x_{s,j}^{out}\overline{A}_{s,j}^{out,n} +$$

$$+\Delta t\sum_{i\in J_{s,in}}\left\{a_{i,n_i+1/2}^+\left(\frac{u_{i,n_i+1/2}^- - a_{i,n_i+1/2}^-}{a_{i,n_i+1/2}^+ - a_{i,n_i+1/2}^-}\right)A_{i,n_i+1/2}^- + a_{i,n_i+1/2}^-\left(\frac{a_{i,n_i+1/2}^+ - u_{i,n_i+1/2}^+}{a_{i,n_i+1/2}^+ - a_{i,n_i+1/2}^-}\right)A_{i,n_i+1/2}^+\right\} - \tag{45}$$

$$-\Delta t \sum_{j \in J_{s,out}} \left\{ a^+_{j,1/2} \left( \frac{u^-_{j,1/2} - a^-_{j,1/2}}{a^+_{j,1/2} - a^-_{j,1/2}} \right) A^-_{j,1/2} + a^-_{j,1/2} \left( \frac{a^+_{j,1/2} - u^+_{j,1/2}}{a^+_{j,1/2} - a^-_{j,1/2}} \right) A^+_{j,1/2} \right\}$$

It is clear that left-hand side of (45) will be nonnegative if

$$\Delta x^{in}_{s,i} \overline{A}^{in,n}_{s,i} + \Delta t\, a^-_{i,n_i+1/2} A^+_{i,n_i+1/2} \geq 0 \tag{46}$$

$$\Delta x^{out}_{s,j} \overline{A}^{out,n}_{s,j} - \Delta t\, a^+_{j,1/2} A^-_{j,1/2} \geq 0 \tag{47}$$

The statement of the theorem follows from (44), (46)-(47) and the CFL restriction.

From (41) one can see that the positivity preserving step depends on ratios $A^{\pm}_{j\mp 1/2}/\overline{A}_j$. Due to bottom topography, irregular channel geometry and possibility of a channel drying these ratios can be much more than unity. In this case, the positive preserving time step restrictions are more severe than the CFL restrictions. In order to overcome these restrictions and guarantee the positivity preserving of our scheme, we adopt a draining time step technique from [41], [25], [35]. The basic idea of this approach is to reduce time step locally only for the cell faces at which the inequalities (44), (46)-(47) do not necessarily hold.

Follow to [35] we introduce the draining time step

$$\Delta t^{drain}_{p,k} = \frac{\Delta x_{p,k} \overline{A}^n_{p,k}}{\max(0, H^{(1)}_{p,k+1/2}) + \max(0, -H^{(1)}_{p,k-1/2})} \tag{48}$$

$$\Delta t^{in,drain}_{s,i} = \frac{\Delta x^{in}_{s,i} \overline{A}^{in,n}_{s,i}}{a^-_{i,n_i+1/2} A^+_{i,n_i+1/2}} \left( \frac{a^+_{i,n_i+1/2} - a^-_{i,n_i+1/2}}{a^+_{i,n_i+1/2} - u^+_{i,n_i+1/2}} \right) \tag{49}$$

$$\Delta t^{out,drain}_{s,j} = \frac{\Delta x^{out}_{s,j} \overline{A}^{out,n}_{s,j}}{a^+_{j,1/2} A^-_{j,1/2}} \left( \frac{a^+_{j,1/2} - a^-_{j,1/2}}{u^-_{j,1/2} - a^-_{j,1/2}} \right) \tag{50}$$

which describes the time when the water contained in cell in the beginning of the time step outflows from it. Now the evolution step in (43), (45) is replaced by

$$\overline{A}^{n+1}_{p,k} = \overline{A}^n_{p,k} - \frac{\Delta t_{p,k+1/2} H^{(1)}_{p,k+1/2} - \Delta t_{p,k-1/2} H^{(1)}_{p,k-1/2}}{\Delta x_{p,k}} \tag{51}$$

$$\sum_{i \in J_{s,in}} \Delta x^{in}_{s,i} \overline{A}^{in,n+1}_{s,i} + \sum_{j \in J_{s,out}} \Delta x^{out}_{s,j} \overline{A}^{out,n+1}_{s,j} = \sum_{i \in J_{s,in}} \Delta x^{in}_{s,i} \overline{A}^{in,n}_{s,i} + \sum_{j \in J_{s,out}} \Delta x^{out}_{s,j} \overline{A}^{out,n}_{s,j} - \tag{52}$$

$$- \sum_{j \in J_{s,out}} \Delta t_{j,1/2} H^{(1)}_{j,1/2} + \sum_{i \in J_{s,in}} \Delta t_{i,n_i+1/2} H^{(1)}_{i,n_i+1/2}$$

The outflow of a cell is the inflow of its neighbor, thus to keep the conservativity of the scheme, we define $\Delta t_{p,k+1/2}$ as

$$\Delta t_{p,k+1/2} = \min(\Delta t, \Delta t^{drain}_{p,m}), \qquad m = k + \frac{1}{2} - \frac{\text{sgn}(H^{(1)}_{p,k+1/2})}{2} \tag{53}$$

where $\Delta t$ satisfying CFL conditions.

For consistency the momentum conservation equation (32) should be rewritten in form

$$\bar{Q}^{n+1}_{p,k} = \frac{\bar{Q}^n_{p,k} - \frac{\Delta t_{p,k+1/2}}{\Delta x_{p,k}} H^{(2,a)}_{p,k+1/2} + \frac{\Delta t_{p,k-1/2}}{\Delta x_{p,k}} H^{(2,a)}_{p,k-1/2}}{1 + \Delta t\, G^{(2)}_{p,k}} - \tag{54}$$

$$- \frac{\Delta t}{\Delta x_{p,k}} \frac{H^{(2,g)}_{p,k+1/2} - H^{(2,g)}_{p,k-1/2} - g\, I_2(x_{p,k-1/2}, x_{p,k+1/2}) + g B_x(x_{p,k-1/2}, x_{p,k+1/2})}{1 + \Delta t\, G^{(2)}_{p,k}}$$

where $H^{(2,a)}_{p,k\pm 1/2}$ and $H^{(2,g)}_{p,k\pm 1/2}$ are the advective and gravity driven parts of the $H^{(2)}_{p,k\pm 1/2}$ which are defined as

$$H^{(2,a)}_{p,k+1/2} = \frac{a^+_{p,k+1/2}}{a^+_{p,k+1/2} - a^-_{p,k+1/2}} \frac{(Q^-_{p,k+1/2})^2}{A^-_{p,k+1/2}} - \frac{a^-_{p,k+1/2}}{a^+_{p,k+1/2} - a^-_{p,k+1/2}} \frac{(Q^+_{p,k+1/2})^2}{A^+_{p,k+1/2}} - \tag{55}$$

$$+ \frac{a^+_{p,k+1/2} a^-_{p,k+1/2}}{a^+_{p,k+1/2} - a^-_{p,k+1/2}} \left[ Q^+_{p,k+1/2} - Q^-_{p,k+1/2} \right]$$

$$H^{(2,g)}_{p,k+1/2} = g\, \frac{a^+_{p,k+1/2} I_1(x_{p,k+1/2}, h^-_{p,k+1/2}) - a^-_{p,k+1/2} I_1(x_{p,k+1/2}, h^+_{p,k+1/2})}{a^+_{p,k+1/2} - a^-_{p,k+1/2}} \tag{56}$$

### 4.8. Well-balancing

The system (1)-(2) admits smooth steady-state solutions, satisfying

$$Q = uA = Const, \qquad E = \frac{1}{2}\frac{Q^2}{A^2} + g(h + B) = Const$$

as well as nonsmooth steady-state solutions. One of the most important steady state solutions is a trivial stationary one (lake at rest)

$$Q = 0, \qquad h + B = Const$$

A numerical scheme that exactly preserves steady flow at rest is called well-balanced.

**Theorem 2.** *Consider the semi-discrete central-upwind scheme (51)-(52), (54) with the piecewise linear reconstruction described in section 4.3 and the discretization of the source terms (14)-(16). Assume that the numerical solution $U(t^n)$ corresponds to the steady state at rest, then $U(t^{n+1})=U(t^n)$, that is, the scheme is well balanced.*

*Proof.* In result of the reconstruction, we have $w_{j\pm 1/2}^{\mp} = w_j$ and $Q_{j\pm 1/2}^{\mp} = 0$. It is clear that the equations of the mass conservation (51)-(52) are satisfied. Thus, we have to show that in the momentum conservation equations (54), there have to be a balance between the flux gradient and the source term. For the second component $H_{j+1/2}^{(2)}$ from (55)-(56) we obtain

$$H_{j+1/2}^{(2,a)} + H_{j+1/2}^{(2,g)} = gI_1(x_{j+1/2}, h_{j+1/2}) \tag{57}$$

as $h_{j+1/2}^{+} = h_{j+1/2}^{-}$ and $A_{j+1/2}^{+} = A_{j+1/2}^{-}$. On the other hand, from (14)-(16) we have

$$I_1(x_{j+1/2}, h_{j+1/2}) - I_1(x_{j-1/2}, h_{j-1/2}) - I_2(x_{j-1/2}, h_{j-1/2}) + B_x(x_{j-1/2}, h_{j-1/2}) =$$

$$= \int_{x_{j-1/2}}^{x_{j+1/2}} \frac{d}{dx} \int_0^h (w - B - y)\sigma(x, y) dy dx - \int_{x_{j-1/2}}^{x_{j+1/2}} \int_0^h (w - B - y)\frac{d}{dx}\sigma(x, y) dy dx + \tag{58}$$

$$+ \int_{x_{j-1/2}}^{x_{j+1/2}} \frac{dB}{dx} \int_0^h \sigma(x, y) dy dx = 0$$

which ensures a balance of source terms with the flux terms in steady states at rest.

## 5. Numerical results

In this section, we examine the scheme accuracy on the problem with smooth solutions and compare numerical results with analytical and measured data as also with the results obtained earlier using some other schemes.

In all tests we assume that water flow is frictionless, the CFL number is 0.5, the acceleration of gravity $g=9.81$, and $\Delta x_j = 0.005$, unless otherwise mentioned.

### 5.1. Accuracy of the scheme in smooth regions

We expect that the scheme (51)-(52), (54) have the first-order accuracy in time and second order in space. We consider test proposed in [22] to check the second order of the scheme in space for smooth flows by evaluation of $L^1$ error over the successively thinner grids.

Let a flat channel has a trapezoidal geometry

$$\sigma(x, y) = 1 + 0.3y$$

Initial water surface and velocity are described by

$$w_0(x) = 1.6 + 0.1\cos\left(\frac{\pi(x - 0.4)}{0.2}\right), \quad u_0(x) = 1, \quad 0 \leq x \leq 1$$

The outflow boundary conditions are set at the channel ends.

We compute $L^1$ error for $\Delta x = 1/N$, $N=80, 160, 320, 640, 1280$ and $2560$ at the final time $T=0.05$. The reference solution was computed with $N=5120$. In order to suppress the scheme error in time, we provide all calculations with time step $\Delta t = 1 \times 10^{-8}$ which is less than $(\Delta x)^2 = (1/5120)^2 \approx 0.38 \times 10^{-7}$. Water surface elevation at different times and $L^1$ error are given in Fig.6, and Table 1, respectively. The results indicate the second-order accuracy in space of the scheme.

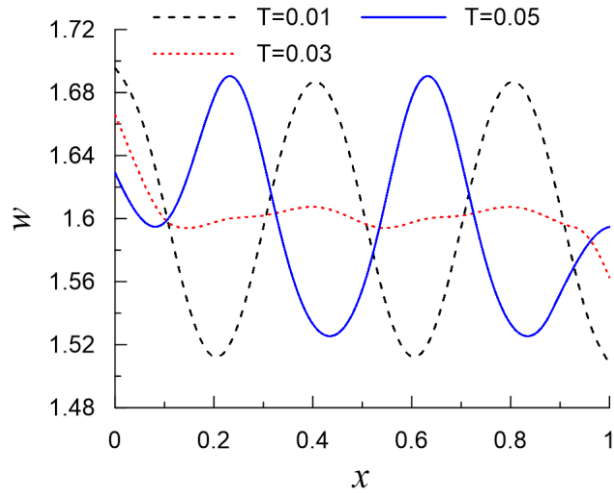

**Fig.6.** Water surface elevation at time t=0.01, 0.03, 0.05 for a trapezoidal channel

| N | w | | Q | |
|---|---|---|---|---|
| | $L^1$-error | order | $L^1$-error | order |
| 80 | 9.60751e-4 | | 1.09671e-2 | |
| 160 | 2.37650e-4 | 2.015 | 2.85182e-3 | 1.943 |
| 320 | 6.19365e-5 | 1.940 | 7.33061e-4 | 1.960 |
| 640 | 1.64387e-5 | 1.914 | 1.89283e-4 | 1.953 |
| 1280 | 4.45586-6 | 1.883 | 4.97988e-5 | 1.926 |
| 2560 | 1.07261-6 | 2.055 | 1.14461e-5 | 2.121 |

**Table 1.** Accuracy test for the scheme in smooth regions. The errors at $T=0.05$.

### 5.2. Large perturbation of rest

The following test is also taken from [22], but with another channel geometry. In this test, propagation of perturbation of a steady state solution at the left part of simulated area is considered. When the perturbation propagates, part of it reflects back and part is transmitted, inundating the right side and leaving the area through the right boundary. According to [22], the bottom topography is described by a piecewise spline. At the beginning we construct the cubic spline of defect one with knots (0,0.3), (0.05,0.3), (0.1,0.2), (0.15,0.5), (0.3,0.4), (0.4,0.6), (0.75,0.6) and zero-values of the second derivatives at the boundaries. Then the spline is rescaled by a factor of 0.5 for $x \leq 0.53$ and of 1.3 for $x \geq 0.53$. Finally, we assume that near the boundaries the bottom is flat and equals to 0.12 for $x \leq 0.0693$ and 0.8 for $x \geq 0.7386$.

The geometry is given by (Fig.7)

$$\sigma(x, y) = \left(1 + \frac{3}{4}\cos(\pi x)\right)\left(1 - \frac{y}{2}\right), \quad 0 \le x \le 1, \quad 0 \le y \le 2 \tag{59}$$

We consider the left boundary as an open [42], [43] and outflow condition is set up on the right boundary.

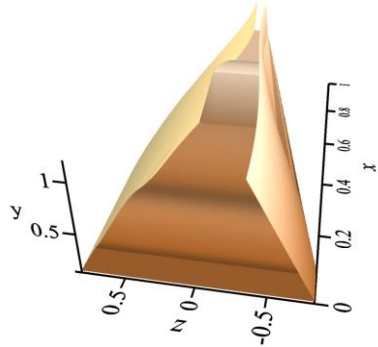

**Fig.7.** 3D view of the channel geometry.

The initial water surface elevation at rest is $w = 0.8$ and a perturbation of size $\varepsilon = 0.3$ is applied on the interval [0.1,0.15]. Propagation of the perturbation at different times are shown in Fig.8. The wave partially reflects back and partially transmits through the bottom ridge, overtopping it and then propagating through the shore. These computed results are similar to results obtained in [22].

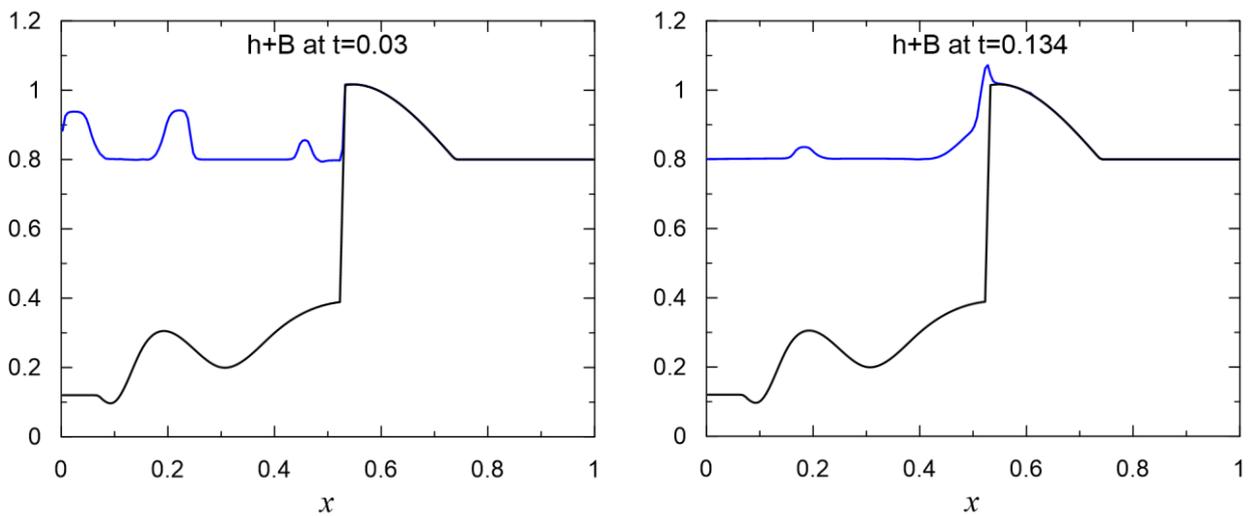

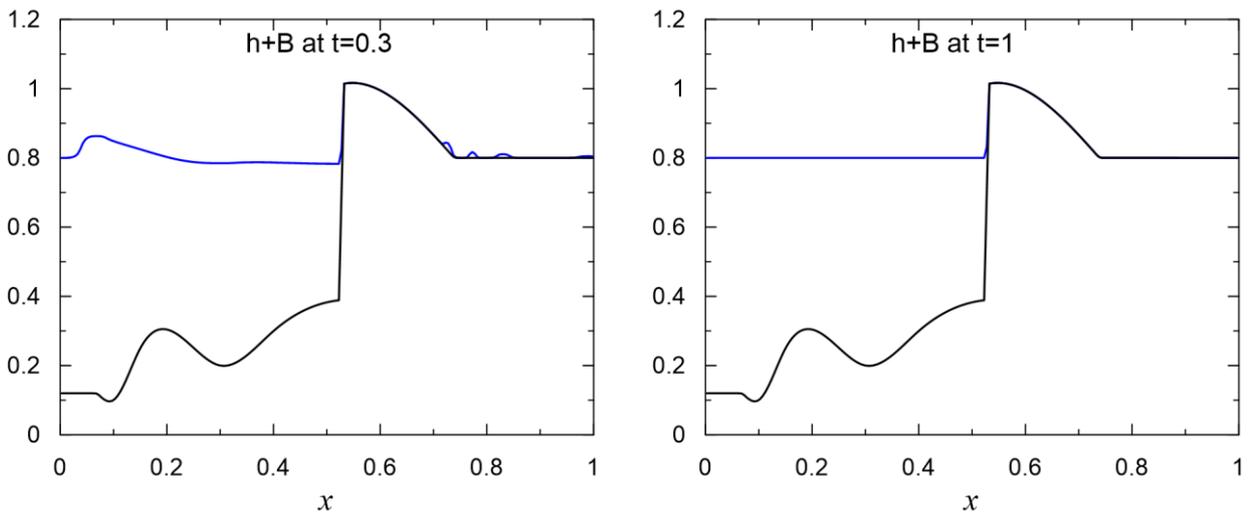

**Fig.8.** Propagation of the perturbation at different times for the conditions of the test proposed in [22].

### 5.3. Convergence to a smooth subcritical steady state

In this example from [23] the topography is described by the cubic spline of defect one with knots (0.2,0), (0.3,0.6), (0.4,0.4), (0.5,0.5), (0.6,0.2), (0.7,0) and zero-value of the second derivatives at the boundaries. The channel geometry is given by (59).

The initial data is

$$w(x,0) = 0.8 \text{ and } u(x,0) = 0$$

At the left (inflow) boundary the discharge $Q_{in} = 0.3343$ and at the right (outflow) boundary the water surface elevation $w = 0.8$ are specified. Hernandez-Duenas and Beljadid in [23] calculated that the subcritical steady-flow with given condition has constant energy $E=10.0748$, in our simulations we obtained almost the same value of the energy $E=10.0307$ with accuracy to four decimal places.

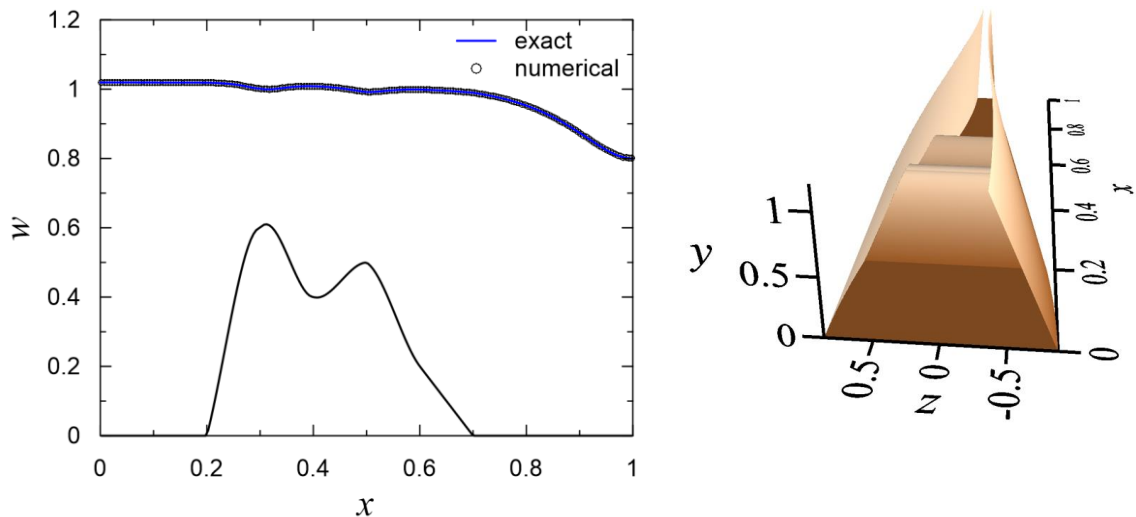

**Fig.9.** Left: Comparison of computed and exact water surface elevation for the subcritical steady state. Right: 3D view of the channel geometry.

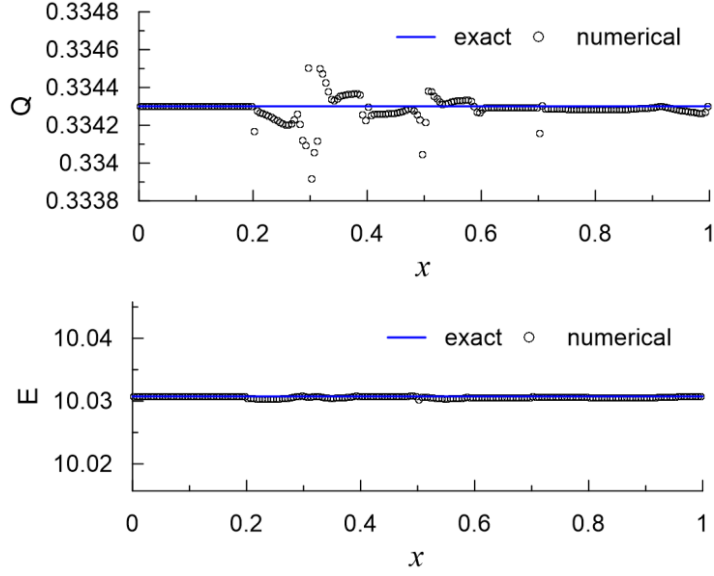

**Fig.10.** Comparison of computed and exact discharges and energies for the smooth subcritical state.

Comparison of computed water surface elevation, discharge, and energy for the smooth subcritical state with their exact values are shown in Fig.9 and Fig.10. In Fig.10 the vertical axis extends on 0.15% of the exact value in each direction. We also as in [23] calculated the maximal error between computed and exact solutions, and relative error in the $L^2$ norm defined as

$$err = \sqrt{\frac{1}{b-a}\int_a^b \left(\frac{f(x)-f_{exact}(x)}{f_{exact}(x)}\right)^2 dx}, \qquad x \in [a,b] \tag{60}$$

We obtained that for the discharge the maximal error is $3.82 \times 10^{-4}$ and the relative error is $1.84 \times 10^{-4}$. For the energy the maximal error is $5.67 \times 10^{-4}$ and the relative error is $2.15 \times 10^{-5}$.

### 5.4. Convergence to a transcritical steady state with shock wave

This test is also taken from [23]. The topography is given by

$$B(x) = \begin{cases} \frac{1}{2} + \frac{1}{2}\cos\left(\pi \frac{x-0.5}{4}\right), & x \in [0.5, 0.9] \\ 0, & x \in [0,1] \setminus [0.5, 0.9] \end{cases}$$

The geometry is

$$\sigma(x,y) = \begin{cases} \frac{1}{2} + \left[\frac{3}{8} - \frac{1}{8}\cos\left(\pi \frac{x-0.7}{0.2}\right)\right]\sqrt{y}, & x \in [0.5, 0.9] \\ \frac{1+\sqrt{y}}{2}, & x \in [0,1] \setminus [0.5, 0.9] \end{cases}$$

The boundary conditions are $Q_{in}=2.5561$ and $w_{out}=1.9968$.

Fig.11 shows that a good agreement is obtained between the analytical solution and the computed water surface elevation. A comparison of the computed discharge and energy with the theoretical results is shown in Fig.12.

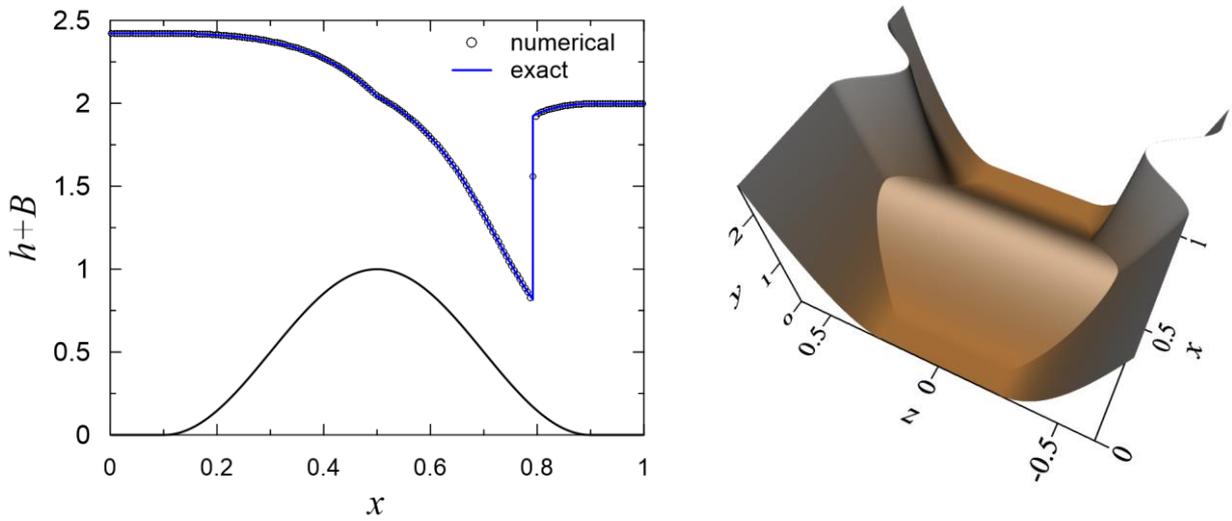

**Fig.11.** Left: Water surface elevation for steady transcritical flow over a bump with a shock. Right: 3D view of the channel geometry.

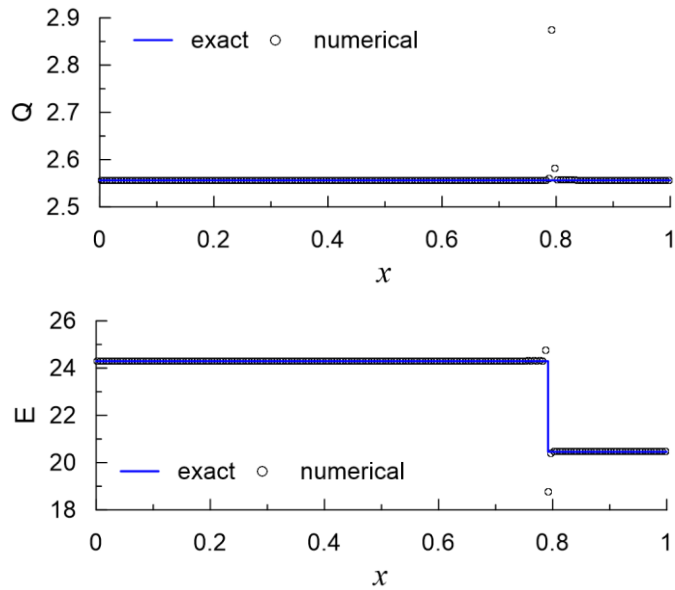

**Fig.12.** Water discharge and energy for steady transcritical flow over a bump with a shock.

*5.5. Dam break problem in a triangular channel*

In this test [44] a dam break problem in a frictionless, horizontal, triangular channel is considered. The channel is 1000 m long with a side slope of 1H:1V. The dam is located in the middle of the channel. Both wet-bed and dry-bed conditions downstream of the dam is simulated. The upstream water depth was 1 m for both cases, and the downstream water depth was 0.1m for the wet-bed case. For the wet-bed case, we compute solution with 400 computational cells. For the dry-bed, we used two computational grids with 400 and 1000 cells.

The exact solutions of dam break problems for a horizontal triangular channel can be found in [45], [46], [47]. Comparison of numerical and exact solution for wet-bed dam break at 80 s after dam removal is shown in Fig. 13. The dry-bed numerical results for water surface and flow rate at 45 s after dam removal are given in Fig.14 and 15. In the dry-bed case, one can see that the most error is at the front of the moving water (Fig.16) which decrease as the cell size is reduced. These results are similar to that obtained by Sanders [48] and Lai and Khan [44] using TVD schemes of second order accurate.

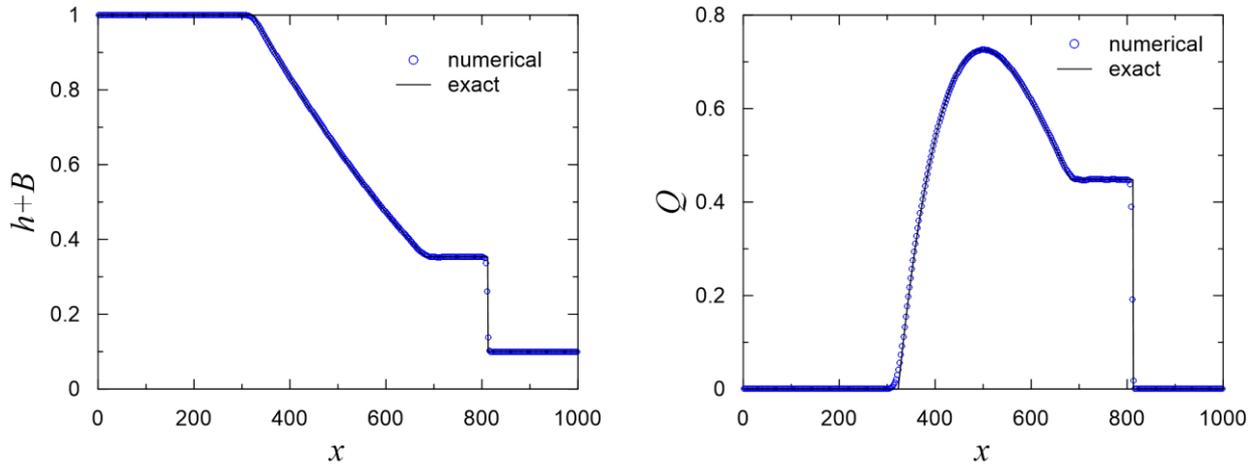

**Fig.13.** Comparison of water surface and flow rate for wet-bed dam break in a triangular channel

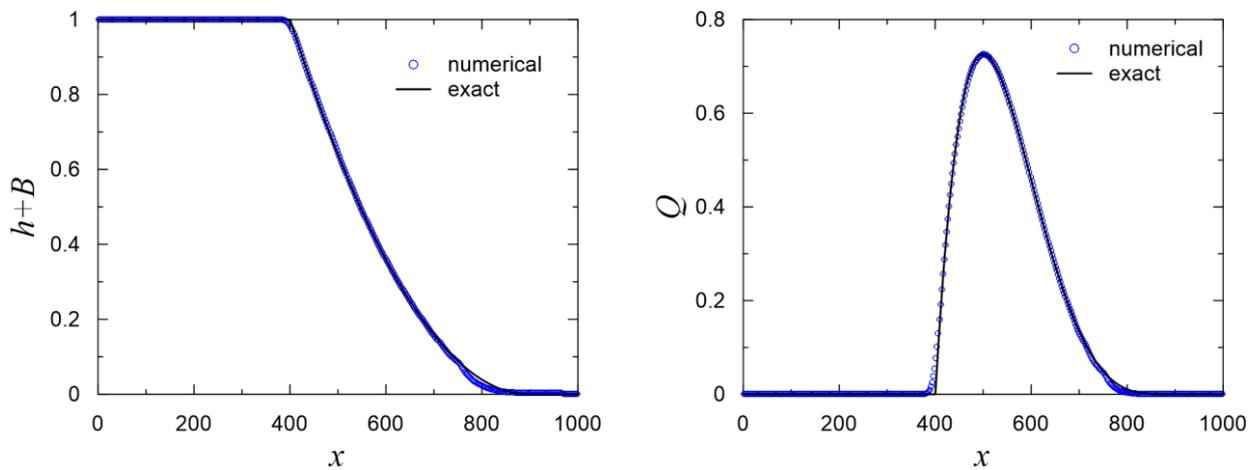

**Fig.14.** Comparison of water surface and flow rate for dry-bed dam break in a triangular channel (N=400)

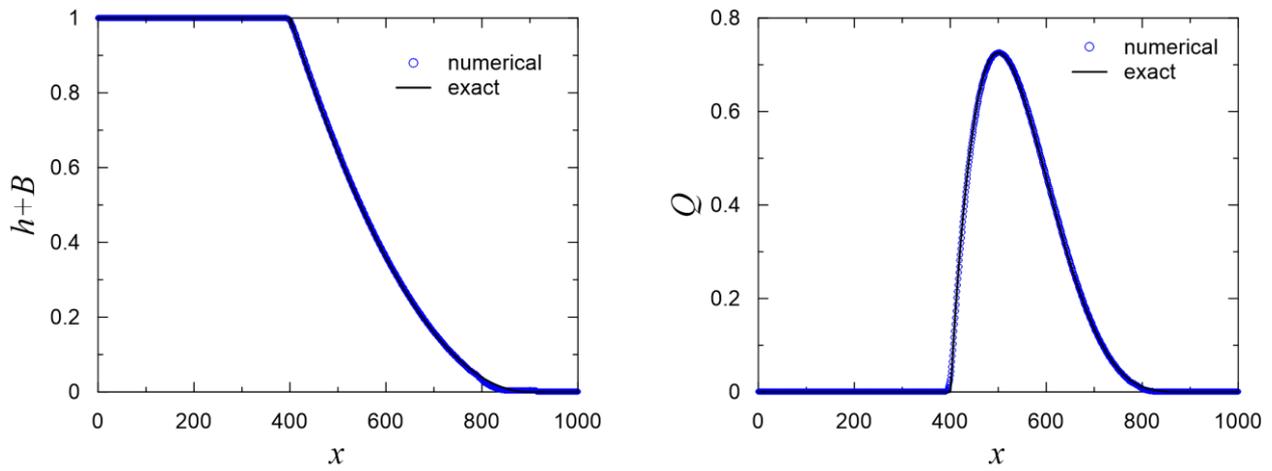

**Fig.15.** Comparison of water surface and flow rate for dry-bed dam break in a triangular channel (N=1000)

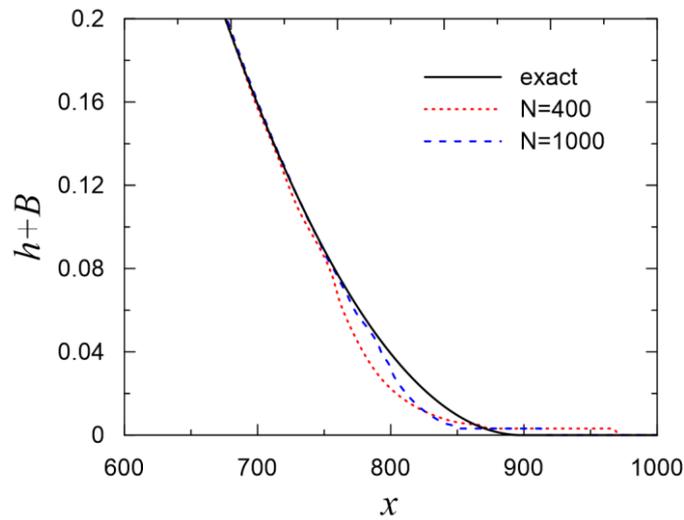

**Fig.16.** Water surface for dry-bed dam break in a triangular channel for different cell sizes

### 5.6. Drain on a non-flat bottom

In this example, taken from [21], [49], a symmetric reservoir is being drained through a parabolic contracting rectangular channel. Due to the symmetry, the flow is computed on half the domain. The contraction is described by the quadratic interpolant through the points (0.25,1.0), (0.5,0.8), and (0.75,1.0), where the first number is *x*-coordinate, and the second one is the width of the channel. The bottom topography consists of one hump (Fig.17)

$$B(x) = \begin{cases} 0.25[1+\cos(\pi(x-0.5)/0.1)], & \text{if } |x-0.5| < 0.1 \\ 0, & \text{otherwise} \end{cases}$$

The left boundary is assumed a line of symmetry of the domain and wall boundary conditions are applied. The right boundary condition is an outflow condition on a dry bed. This boundary condition on the right side of the domain, allows the water that was at rest freely flows through the right boundary into the initially dry region.

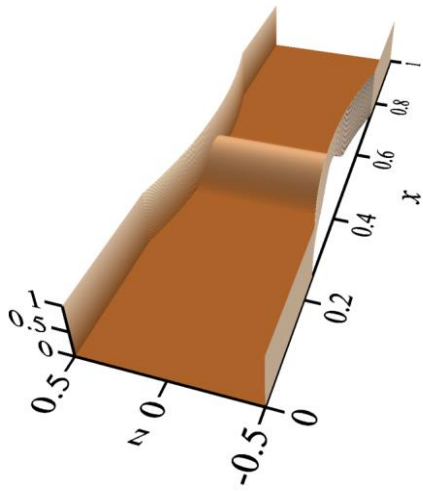

**Fig.17.** 3D-view showing the topography at the bottom and lateral channel walls.

In Fig.18 we show the solution obtained for the initial water surface elevation at 0.8.

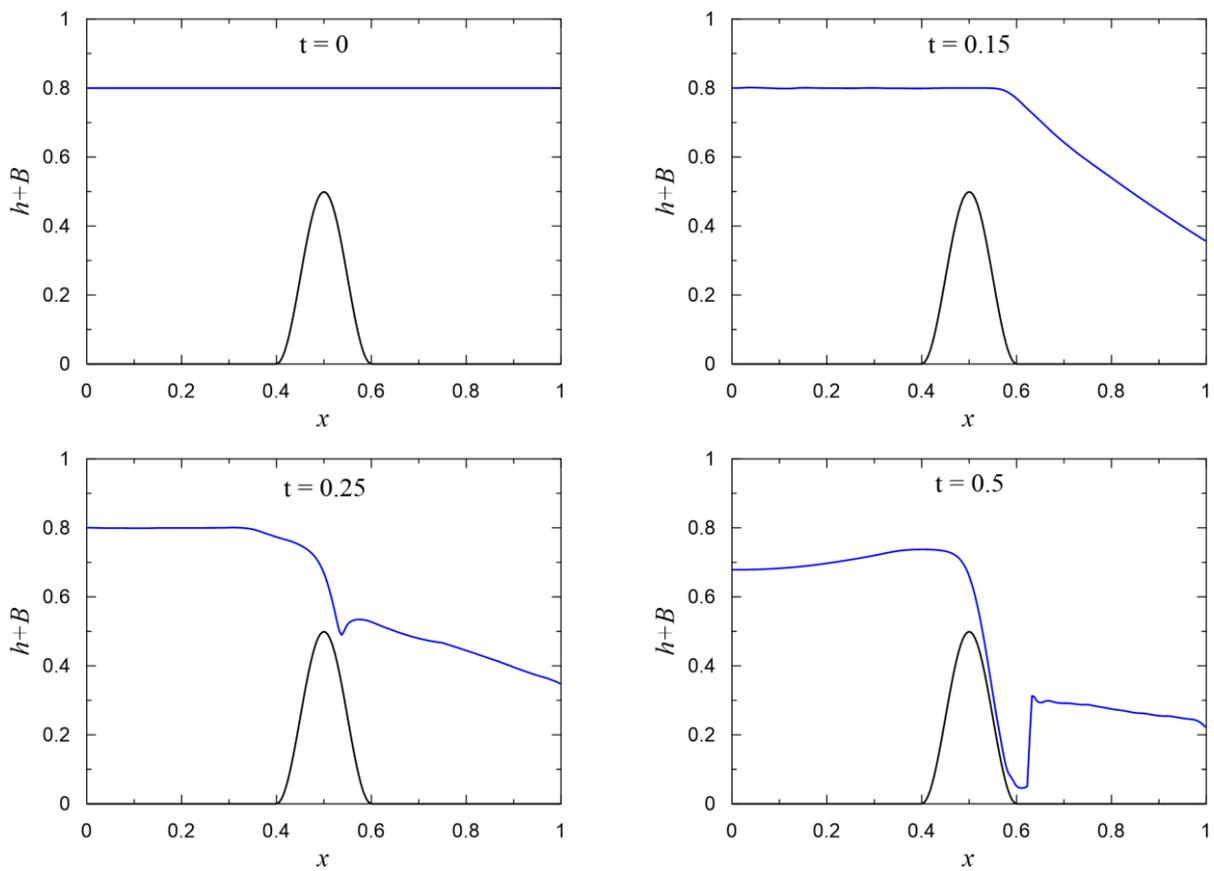

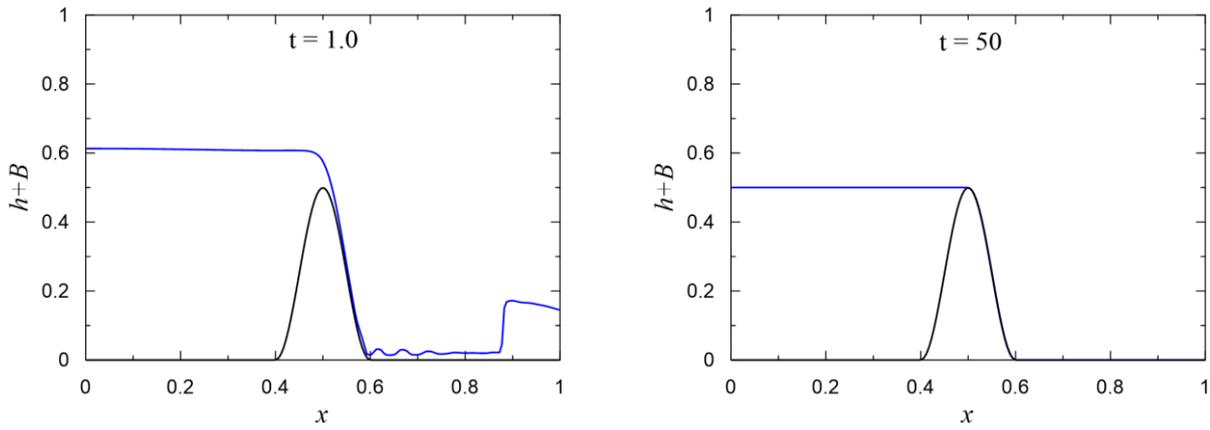

**Fig.18.** Water surface elevation under drain on a non-flat bottom at different times for the test [21, 49].

After drainage beginning the solution converges to a steady-state solution in which water exists only on left hand side of the hump. The obtained simulated results are in good agreement with numerical result in [21].

*5.7. Subcritical dam-break flow in a rectangular channel with a junction*

We consider a subcritical dam-break flow in a frictionless, horizontal, rectangular channel with 34 m long and 3 m wide. The channel is divided into two parts by a dam located 15 m from the left end. Initially, the water at rest to the left from the dam was at a level of 0.5 m and 0.1 m to the right. Vertical walls are placed at the ends of the channel and their height is enough to prevent any spilling of water out of the channel (Fig. 19). The channel cross-sections for numerical solutions comparison G1 and G2 were respectively located 4.4 m and 5.9 m downstream from the dam.

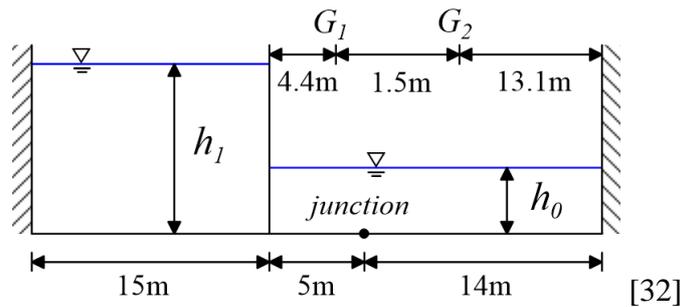

**Fig.19.** Outline of dam break flows in a rectangular channel with a junction

After abrupt dam failure, the water flows downstream until it reaches the vertical wall at the right channel end and reflects from it. The left vertical wall reflects the flow too. The simulations are carried out for a uniform $\Delta x = 0.1$ m and $\Delta t = 0.01$ s. Numerical simulations are performed for a computational grid with and without channel junction that is located 5 m downstream from the dam. We apply two types of junction treatment for a control volume around the channel junction: 1) based on mass conservation (CUJ1) and 2) based on the Saint Venant equations (CUJ2). The numerical solution on the computational grid without junction (CU) we consider as the reference solution.

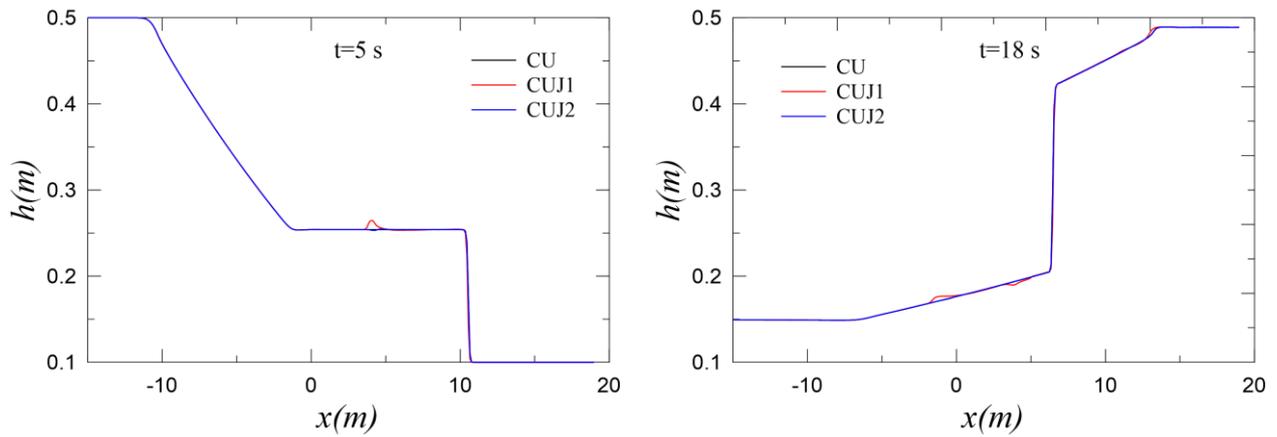

**Fig.20.** Simulated water surface profiles of the entire channel at t = 5 s and 18 s

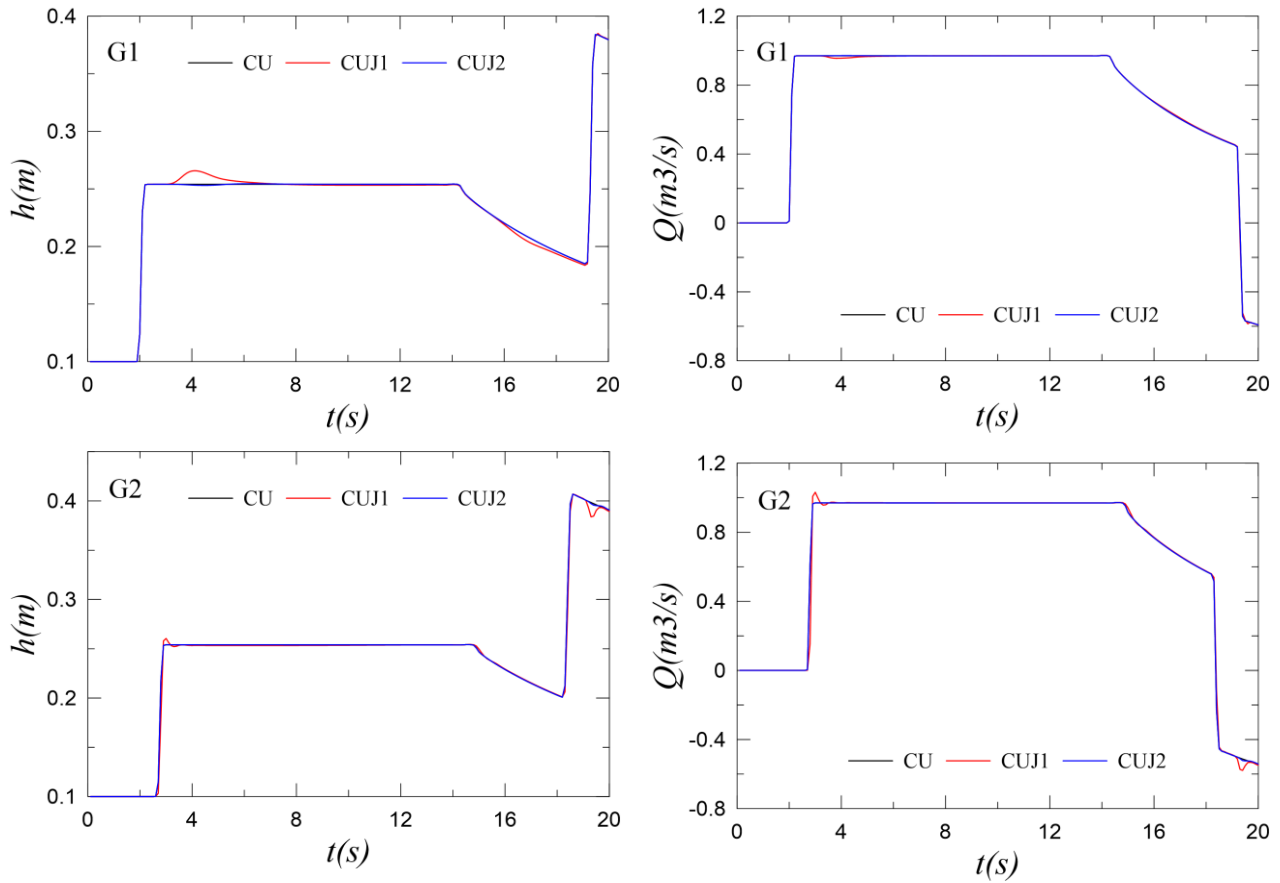

**Fig.21.** Comparisons of the simulated water depths and discharges with different junction treatments at selected cross-sections.

Simulated water profiles for two times at 5 s and 18 s are given in Fig.20. For different junction treatments, comparisons of the simulated water depths at two channel cross-sections G1 and G2 are shown in Fig.21. Applying only the mass conservation equation for a channel junction can lead to small errors in the modeling of a subcritical open water flows.

### 5.8. Dam-break flow in a rectangular channel with a loop

Purpose of this test is to compare different treatments of an open channel junction for subcritical and supercritical flows. Similarly, to the previous section we consider a dam break problem in a frictionless, horizontal, rectangular channel with a loop (Fig.22). There are two vertical walls at ends of the channel located at *x*=-15 m and *x*=19 m. The channel width before and after the loop is 3 m. The width of the lower loop channel (channel B in Fig.22) is 2 m, and the upper (channel C) is 1 m. The channel is divided into two parts by a dam at *x*=0 m. Initially, the upstream water depth is 1 m (supercritical flow), or 0.5 m (subcritical flow), downstream water depth is 0.1 m. The dam is then breached, and instantaneously water discharges from the higher to the lower level as a downstream-directed bore while a depression wave propagates upstream. The flow is reflected when the water front reaches the walls at the ends of the channel. Four channel cross-sections for comparison of numerical solutions G1, G2, G3, and G4 were located at different parts of the channel as shown in Fig.22.

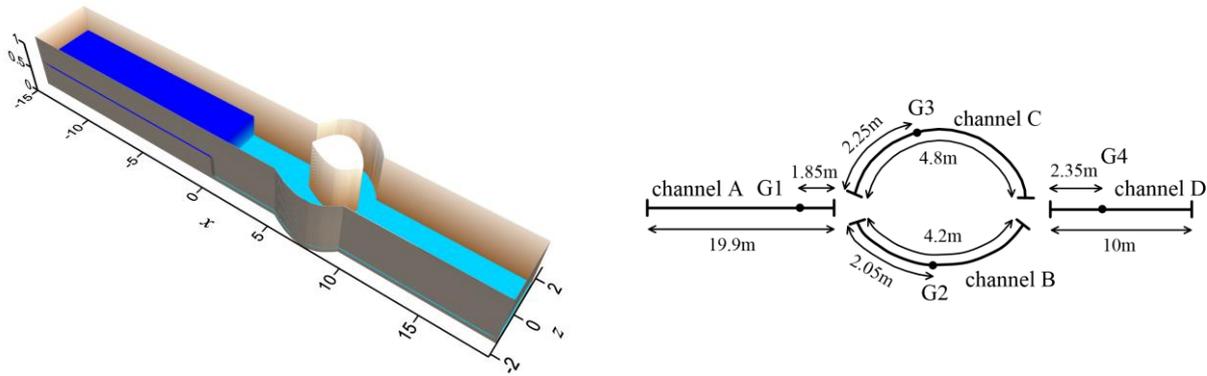

**Fig.22.** Sketch of modeling area for dam break problem in a rectangular channel with a loop.

The simulations are performed with $\Delta x = 0.1$ m and $\Delta t = 0.01$ s. 2D model COASTOX [50], [51], [52] and 1D model are used to simulate dam break problem. COASTOX is an unstructured finite volume model for solving the shallow water equations by using the Godunov scheme. COASTOX numerical solutions are the reference for comparison of the 1D numerical solutions. In 1D simulations, junction treatment in channel network is represented by two kinds of model: 1) mass conservation (CUJ1), and 2) mass and momentum conservation (CUJ2). These two junction models are compared at the four cross-sections with each other and with respect to the averaged over the cross-section numerical solution of the 2D model.

Comparison of the simulated water depths by using two junction models at the four cross-sections for the subcritical flow is present in Fig.23. Relative errors in $L^1$ and $L^2$ norms (60) between 1D and 2D solutions at the cross-sections for the subcritical flow are given in Table 2.

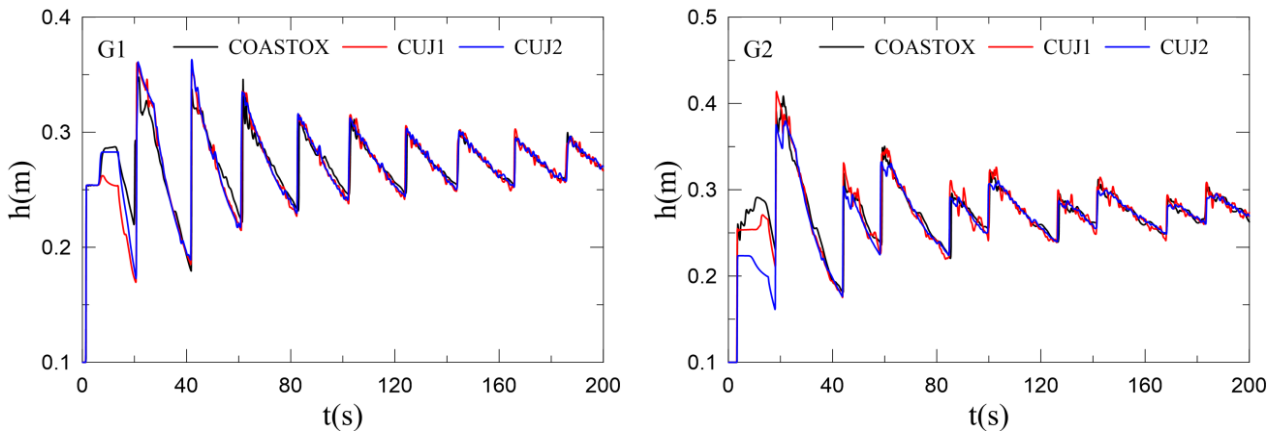

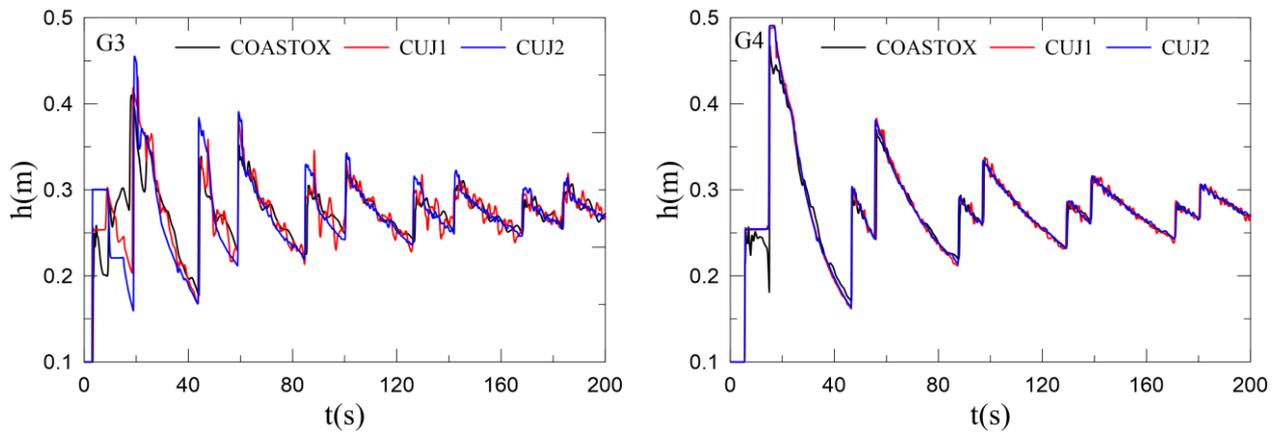

**Fig.23.** Comparison of the simulated water depths at channel cross-sections G1, G2, G3, and G4 for the subcritical flow.

|  | Relative error in $L^1$ norm | | Relative error in $L^2$ norm | |
| --- | --- | --- | --- | --- |
|  | CUJ1 | CUJ2 | CUJ1 | CUJ2 |
| G1 | 0.0298 | 0.0249 | 0.0626 | 0.0578 |
| G2 | 0.0284 | 0.0384 | 0.0506 | 0.0793 |
| G3 | 0.0520 | 0.0725 | 0.0970 | 0.1473 |
| G4 | 0.0230 | 0.0216 | 0.0524 | 0.0711 |

**Table 2.** Relative errors in $L^1$ and $L^2$ norms between 1D and 2D solutions at the cross-sections for the subcritical flow.

CUJ1 produce in similar but the smaller relative errors than CUJ2 for the subcritical flow in this test. The same time the results of CUJ1 much better describe the dynamics of the first wave in the section G2.

For the supercritical flow in a channel network CUJ1 produces not a physical numerical solution. Comparison between CUJ2 and COASTOX numerical solutions for the supercritical flow at the four gages is shown in Fig.24.

Agreement between 1D and 2D numerical results for subcritical and supercritical flows is generally satisfactory when only CUJ2 approach is used for the supercritical flows.

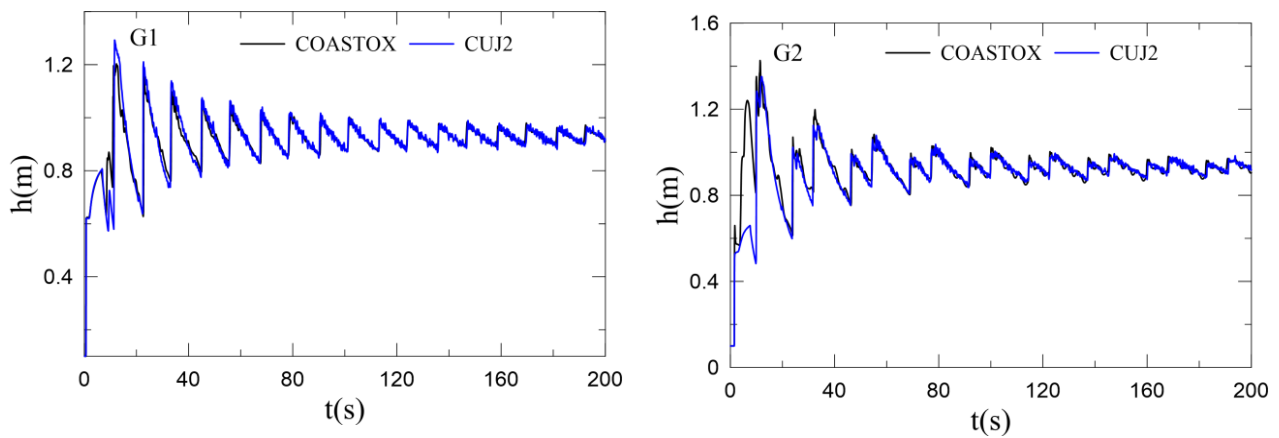

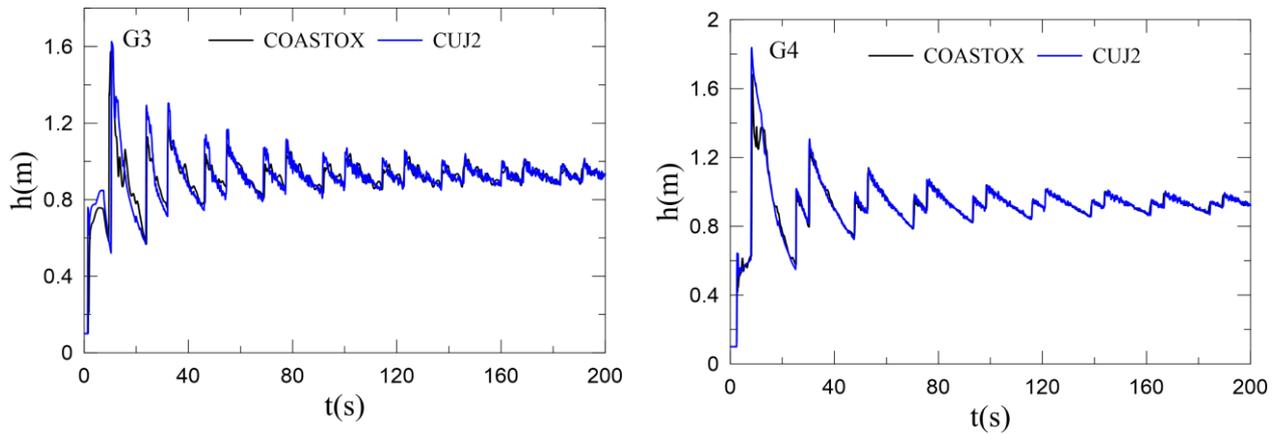

**Fig.24.** Comparison of the simulated water depths by CUJ2 and COASTOX at the channel cross-sections G1, G2, G3, and G4 for the supercritical flow.

As one might suppose, greatest error in the numerical solution of a 1D model in comparison with a 2D model is observed at the bore front under passing of the channel junction. Therefore, it is better to apply 2D or 3D models to simulate a bore propagation in a multiply-connected channel network.

### 5.9. Lower Danube River

As the test of the applicability of the presented above numerical scheme for the large-scale natural open flow we simulate the river network of the Lower Danube from 01/01/2000 to 31/12/2000 [52]. The river network consists of 29 links and 26 nodes (Fig.26, left). Blue and green triangles denote inflow nodes, yellow squares indicate outflow nodes and red circles mark channel junctions. Water gages D1-D10 are located at points indicated by red circles shown in Fig.26 (right). The total length of the main channels is more than 900 km. The total length of the main channels between water gages D1 and D10 is more than 500 km. For this part of the river network, we have only 10 channel cross-sections located at the water gages.

The river network is discretized by a non-uniform grid with cell sizes varied from 1061.37 to 3948.45 m. The computational grid includes 321 cells and 338 cell interfaces. The CFL number is set to 0.8. We will use to simulate water flow in the river network two models: 1) the described in this paper central-upwind scheme with the junction treatment based on continuity equation (CUJ1) and 2) analog of the CHARIMA model [4]. CHARIMA simulates water flow in a channel by applying the Preissmann implicit finite-difference scheme [7]. This scheme is a weighted four-point scheme. At a junction node, it is assumed that: 1) inflow discharges should be equal outflow discharges from all tributaries at the junction, and 2) the water levels at the ends of linked channels are equal at the junction.

The friction slope is evaluated from Manning's equation (3). The hydraulic radius is calculated as A/P, where P is the wetted perimeter. Note that CUJ1 and CHARIMA have different difference approximations of the friction slope.

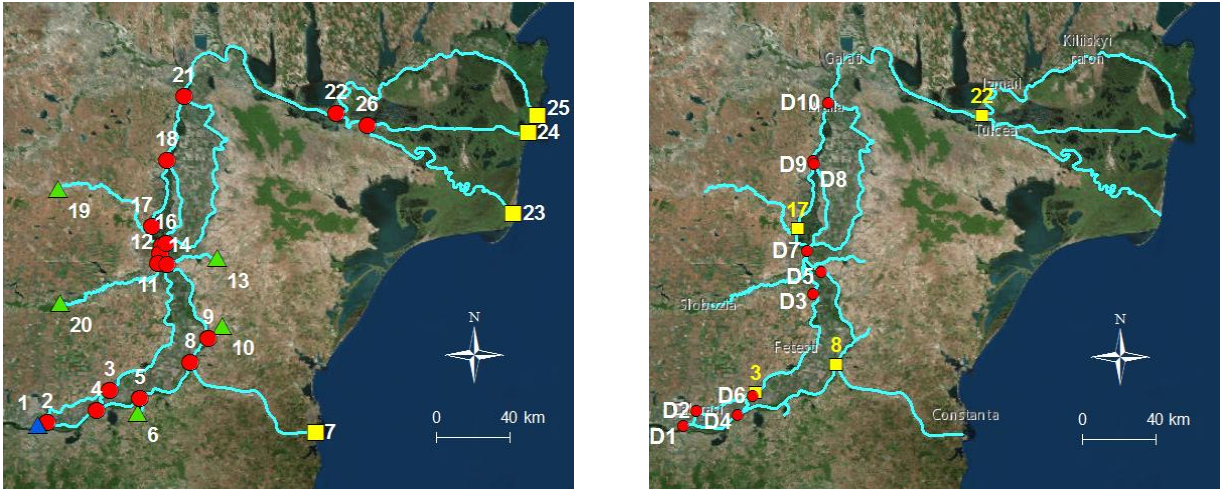

**Fig.25.** River network of the Lower Danube River. Left: The inflow nodes of the river network are indicated by triangles (blue and green) and outflow nodes are indicated by yellow squares. Red circles denote the river channel junctions with their node numbers. Right: Red circles denote location of the water gage stations D1-D10. Yellow squares denote the junction nodes in which the water surface levels calculated by CUJ1 and CHARIMA models are compared.

Water discharge is specified at the inflow nodes. At that, the water discharge at a node indicated by green in Fig.26 (left) is simulated as a lateral inflow to the corresponding channel junction. Water surface elevation is set up at the outflow nodes (yellow squares).

The calibration of CUJ1 is performed by adjusting Manning's coefficient $n$ inside the river network by using the CMA Evolution Strategy [53]. We assumed that $0.01 \leq n \leq 0.04$ and for its evaluation we adopted pCMALib code [54]. The Nash-Sutcliffe model efficiency coefficient (NSE) is applied to assess approximation of the measured water surface elevation at a water gage by the simulated water surface elevation. We took data of water surface level measurements at all water gages D1-D10 for the Manning's coefficient evaluation. Thus, the objective functional is taken as

$$\Im = \sum_{i=1}^{10} NSE_i = \sum_{i=1}^{10}\left[1 - \frac{\sum_{j=1}^{m}\left(w_i(t_j) - w_i^0(t_j)\right)}{\sum_{j=1}^{m}\left(w_i^0(t_j) - \overline{w}_i^0\right)}\right] \to \max$$

where $\overline{w}_i^0$ is the mean of observed water surface elevations at a $i$-th water gage; $w_i^0(t_j)$ is observed water surface elevation at a $i$-th water gage in time $t_j$; $w_i(t_j)$ is modeled water surface elevation at a $i$-th water gage in time $t_j$; m is number of times. When the functional $\Im$ reached a value of 9.5527 we interrupted the calibration. Values of the Nash-Sutcliffe coefficients at gages D1-D10 are given in Table 3.

Postcalibrated values of Manning's n were used in simulations of the two models under same computational grid, initial and boundary conditions. Comparison of simulated and observed water surface elevations at the Borcea (D3), Izvoarele (D4), Gropeni (D8) and Braila (D10) water gages are shown in Fig.26. Comparison of the water surface elevations at junctions 3, 8, 17 and 22 (Fig.25, right) calculated by the CUJ1 and CHARIMA are presented in Fig.27.

|     | D1 | D2 | D3 | D4 | D5 | D6 | D7 | D8 | D9 | D10 |
| --- | --- | --- | --- | --- | --- | --- | --- | --- | --- | --- |
| NSE | 0.9509 | 0.9496 | 0.9799 | 0.9779 | 0.9716 | 0.9592 | 0.8956 | 0.9834 | 0.9884 | 0.8960 |

**Table 3**. Values of the NSE at water gages D1-D10.

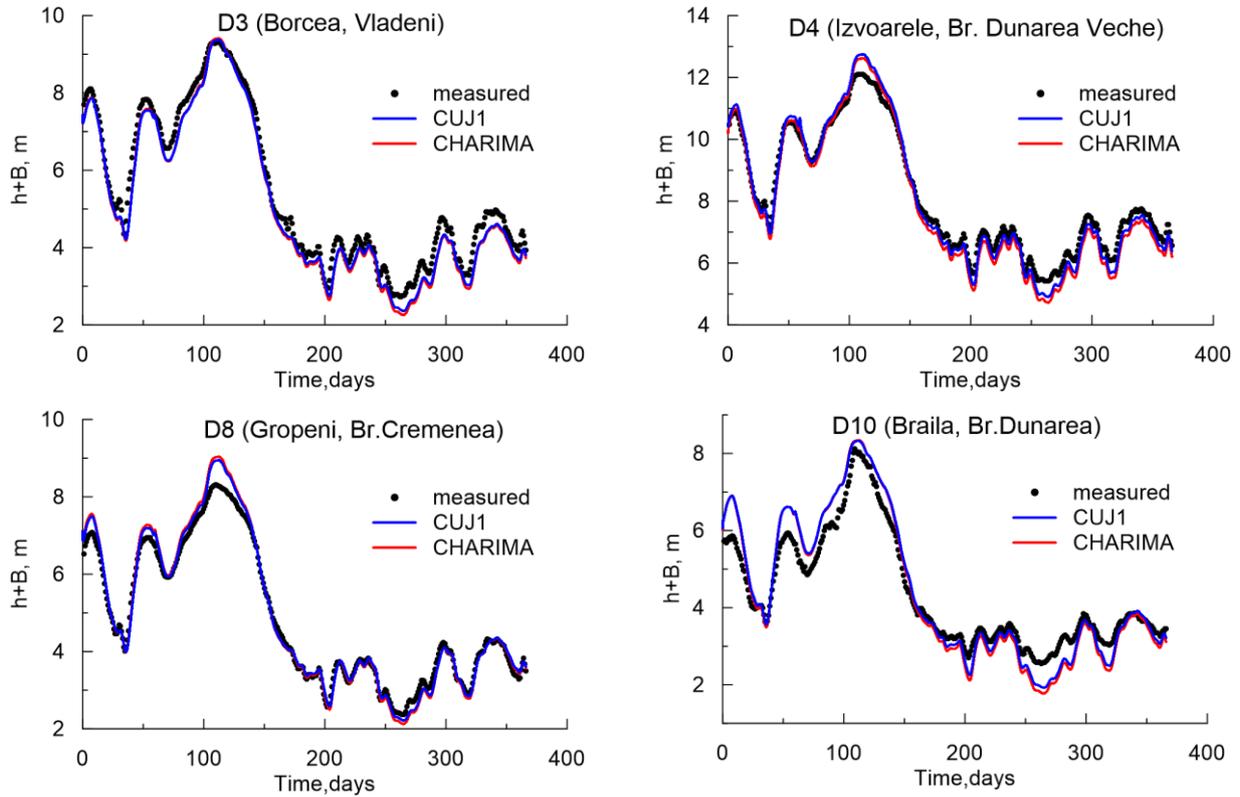

**Fig.26.** Comparison of simulated and observed water surface elevations at various gage stations.

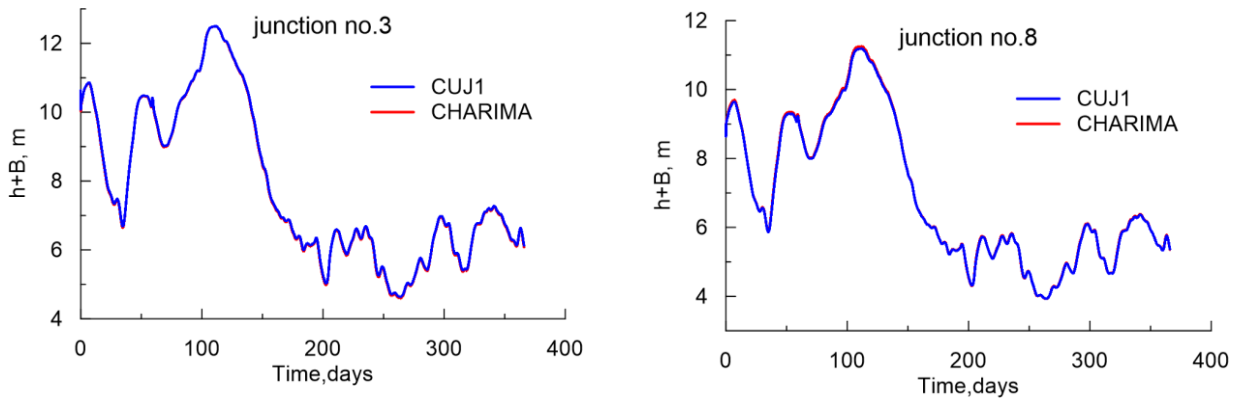

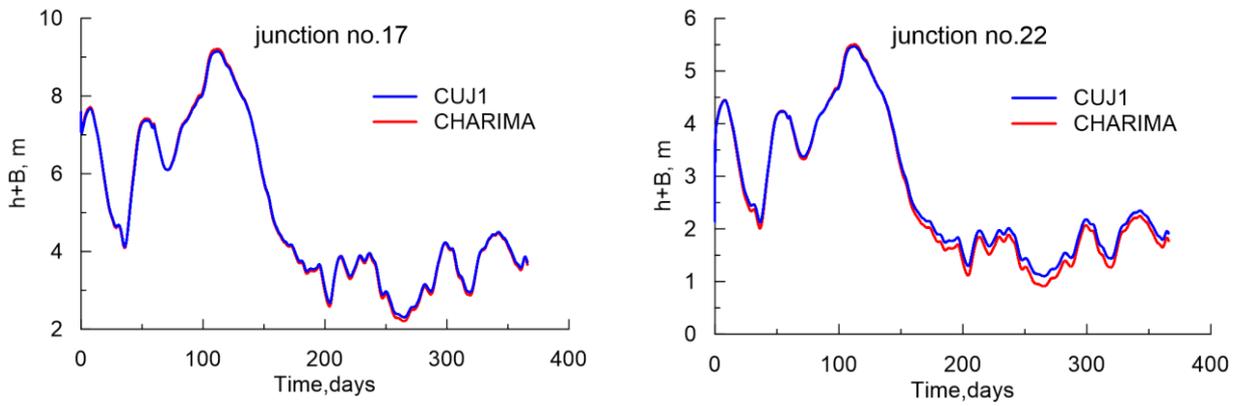

**Fig.27.** Comparison of the water surface elevations at junctions calculated by CUJ1 and CHARIMA.

A good agreement between the water surface elevations at the channel junctions calculated by CUJ1 and CHARIMA models indicates that continuity equation at a junction for a subcritical flow can be used instead a traditional internal boundary condition based on equality of the both discharges and water levels. That provide possibility to apply the described central-upwind scheme CUJ1 for a multiply-connected channel network for a subcritical flows through a junctions.

### *5.10. Inundation of a dry channel network*

In this section we simulate inundation of a dry channel network which consists of eight inclined rectangular channels and four junction nodes (Fig.28). We will denote each channel by the node names bounded this channel. In channel name the left letter will mark the upper end of the channel. Upper ends of channels *AB*, *EF* and *HG* are bounded by vertical walls. Characteristics of network channels are given in Table 4. Note that in the junction *F* channel ends connected with it have different elevations.

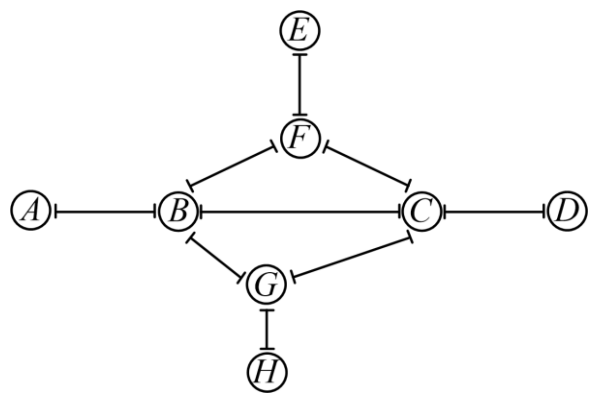

**Fig.28.** Sketch of a channel network.

| Channel | Length, m | Width, m | Elevation of channel ends, m | |
|---|---|---|---|---|
| | | | Upper | Lower |
| AB | 10.0 | 0.1 | 0.25 | 0.20 |
| BC | 20.0 | 0.1 | 0.20 | 0.10 |
| CD | 10.0 | 0.1 | 0.10 | 0.0 |
| EF | 20.0 | 0.2 | 0.35 | 0.20 |
| FB | 15.0 | 0.2 | 0.25 | 0.20 |
| FC | 15.0 | 0.2 | 0.30 | 0.10 |
| HG | 10.0 | 0.2 | 0.35 | 0.30 |
| GB | 10.0 | 0.2 | 0.30 | 0.20 |
| GC | 20.0 | 0.2 | 0.30 | 0.10 |

**Table 4**. Characteristics of network channels.

At the upper ends of the channels *AB*, *EF* and *HG* there are water sources discharges (m³/s) of which are calculated by formulas

$$Q^A(t) = \begin{cases} 0.002\sin(\pi(t-1)/60) & \text{if } 1 \le t \le 61 \\ 0 & \text{otherwise} \end{cases}$$

$$Q^E(t) = \begin{cases} 0.0015\sin(\pi(t-1)/30) & \text{if } 1 \le t \le 31 \\ 0.01\sin(\pi(t-150)/60) & \text{if } 150 \le t \le 210 \\ 0 & \text{otherwise} \end{cases}$$

$$Q^H(t) = \begin{cases} 0.0015\sin(\pi(t-50)/80) & \text{if } 50 \le t \le 130 \\ 0 & \text{otherwise} \end{cases}$$

Such geometry of the channel network and the time distribution of the water sources discharges were selected for this test to generate the consequence of the waves inundating and drying the network.

The channel network was discretized by a uniform grid with cell size of 0.2 m. The Manning's roughness coefficient was equal to 0.01 (s/m$^{1/3}$). The outflow boundary condition was specified at the node *D*. Numerical simulations were carried out for the CFL number equals to 0.9. We used CUJ2 model to simulate open water flow in the channel network under inundation.

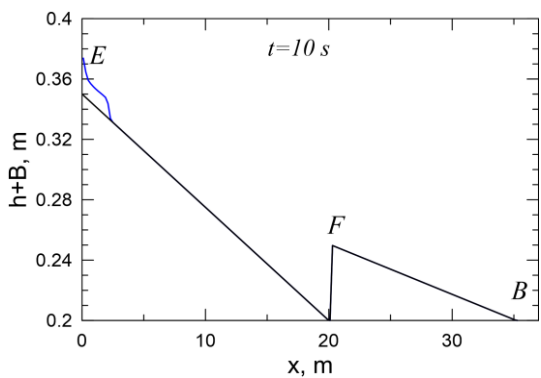
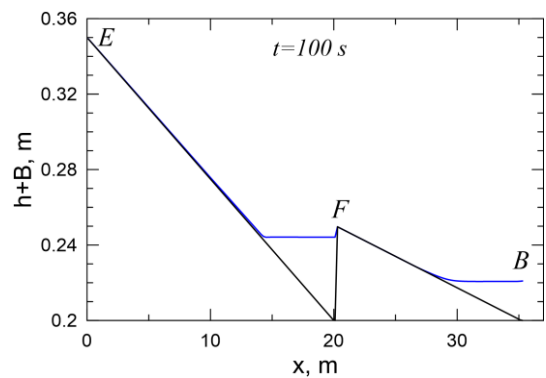

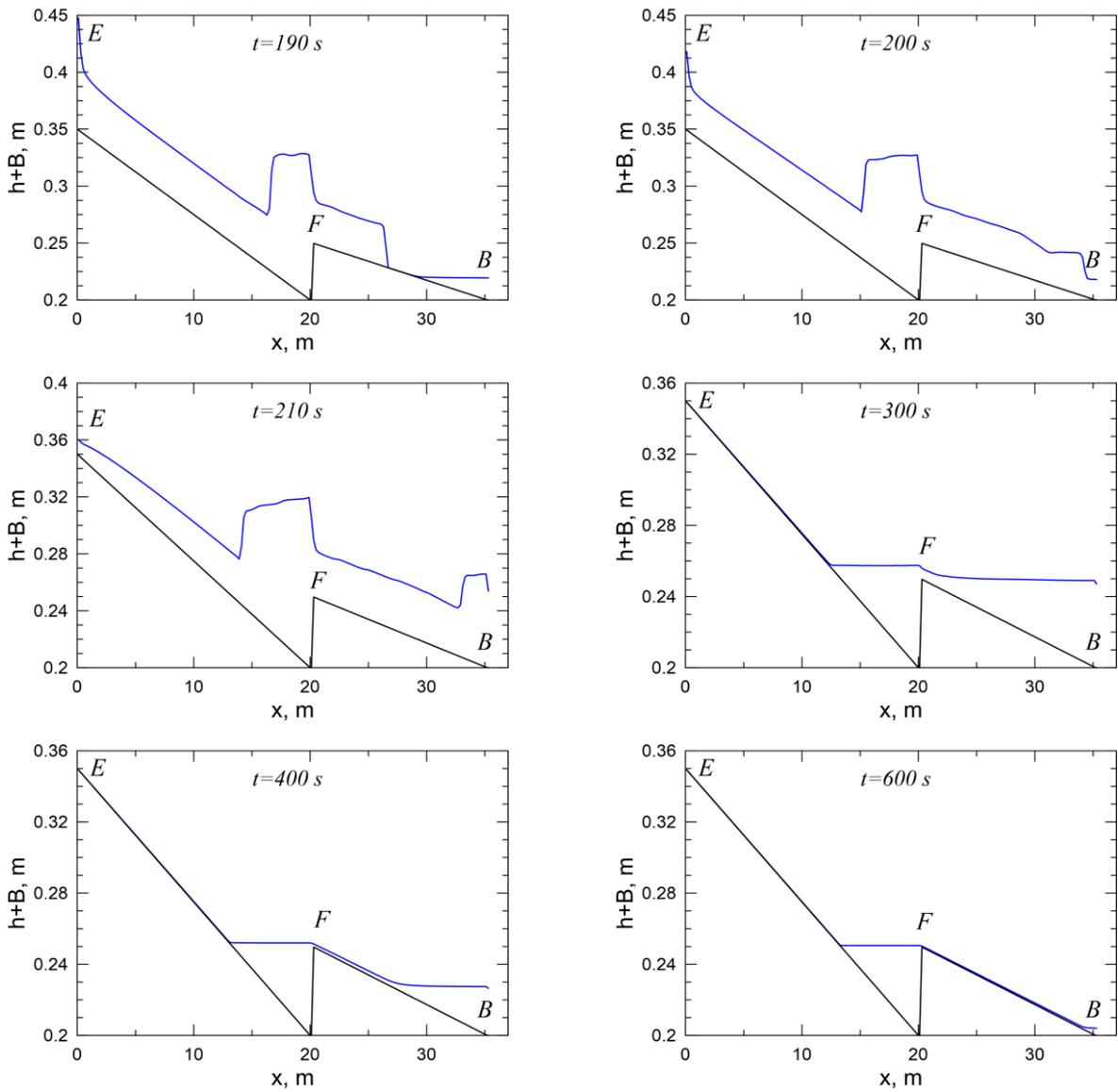

**Fig.29.** Water free surface elevations at different times along channels *EFB*.

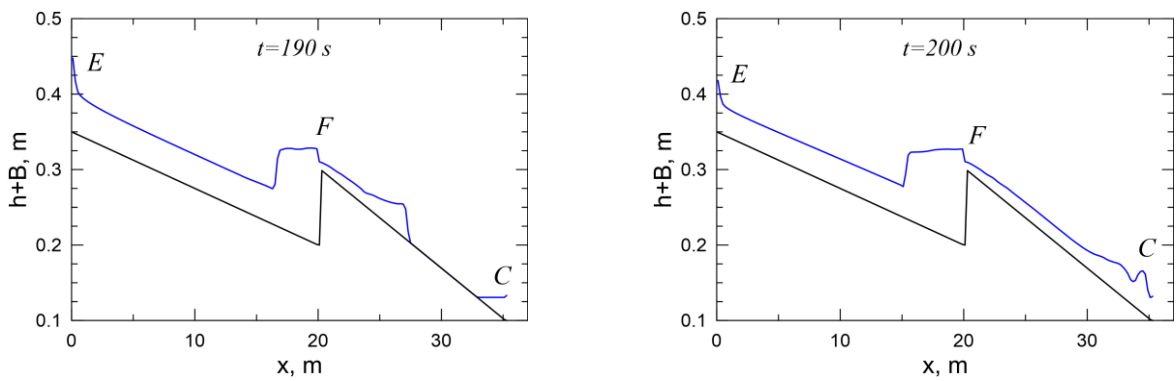

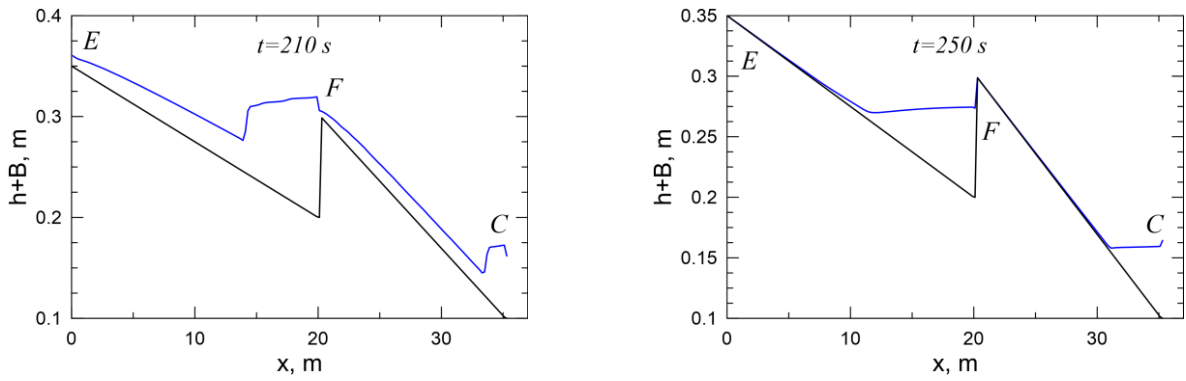

**Fig.30.** Water free surface elevations at different times along channels *EFC*.

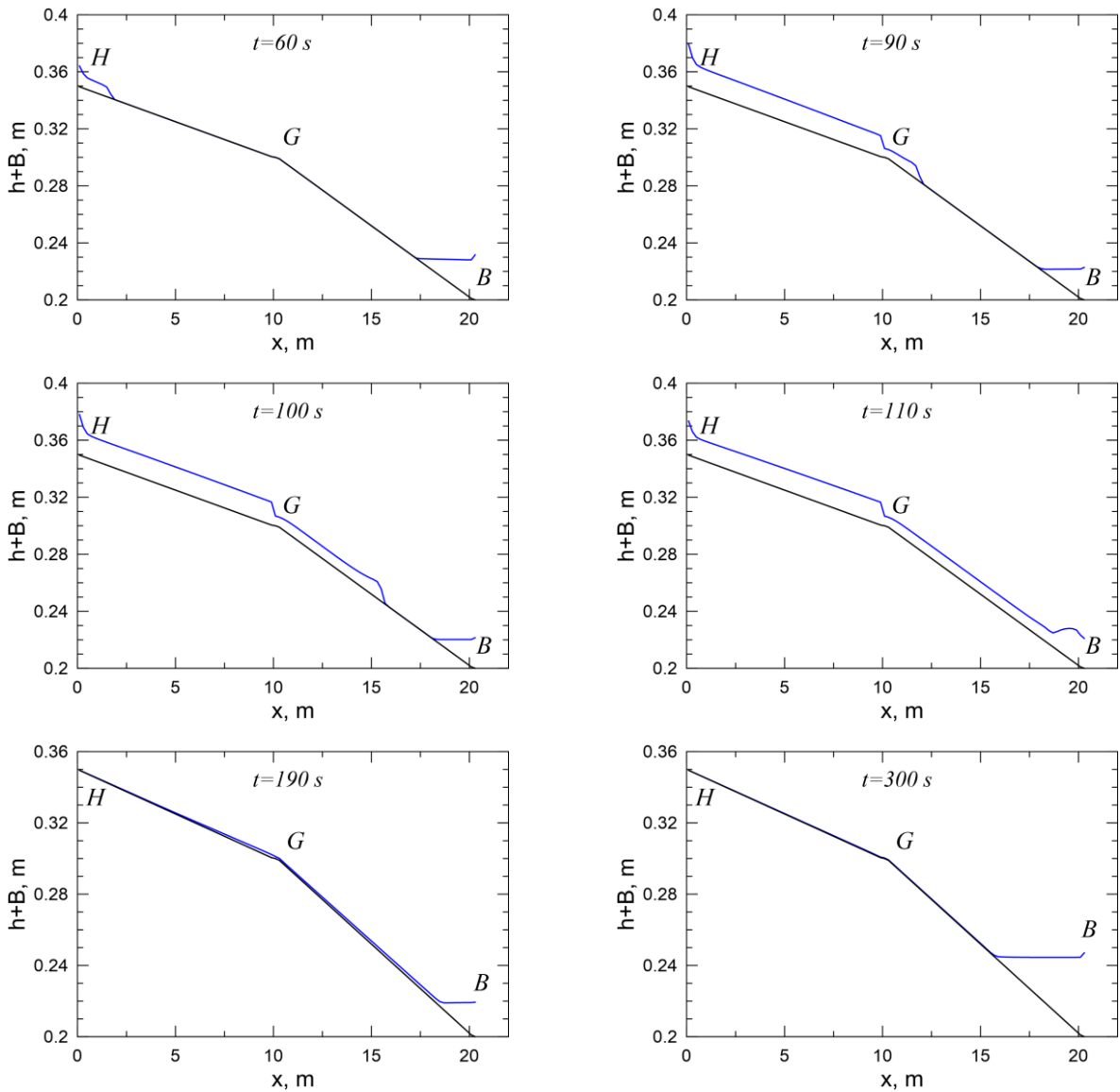

**Fig.31.** Water free surface elevations at different times along channels *HGB*.

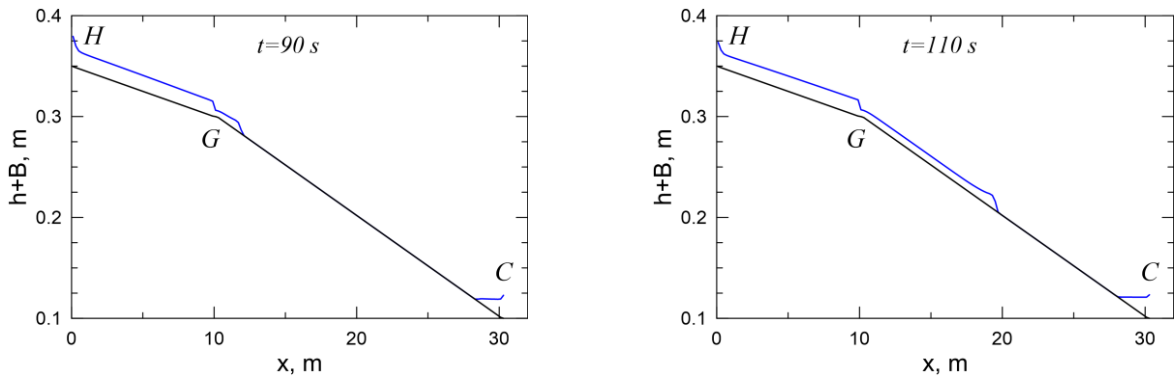

**Fig.32.** Water free surface elevations at different times along channels *HGC*.

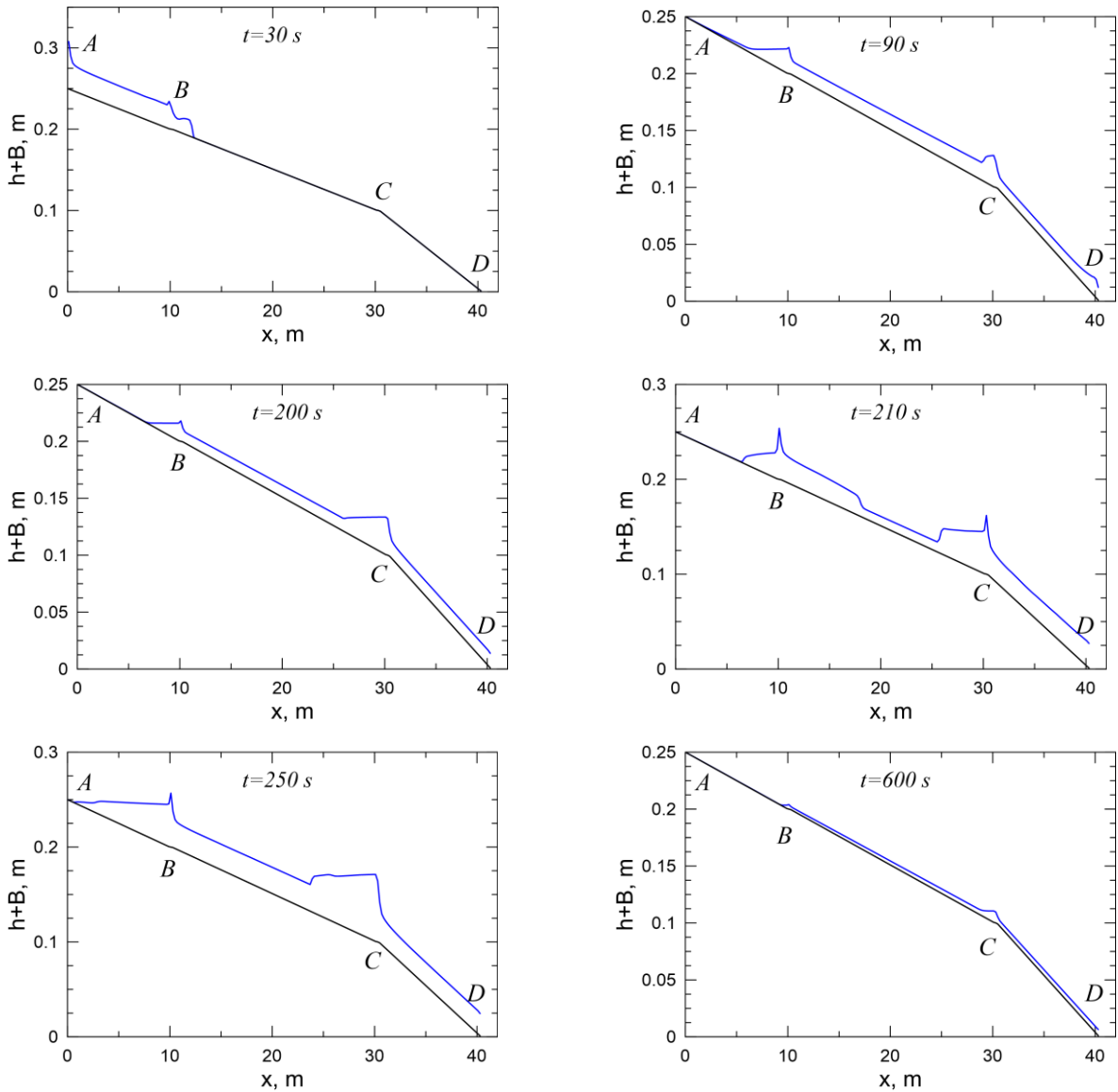

**Fig.33.** Water free surface elevations at different times along channels *ABCD*.

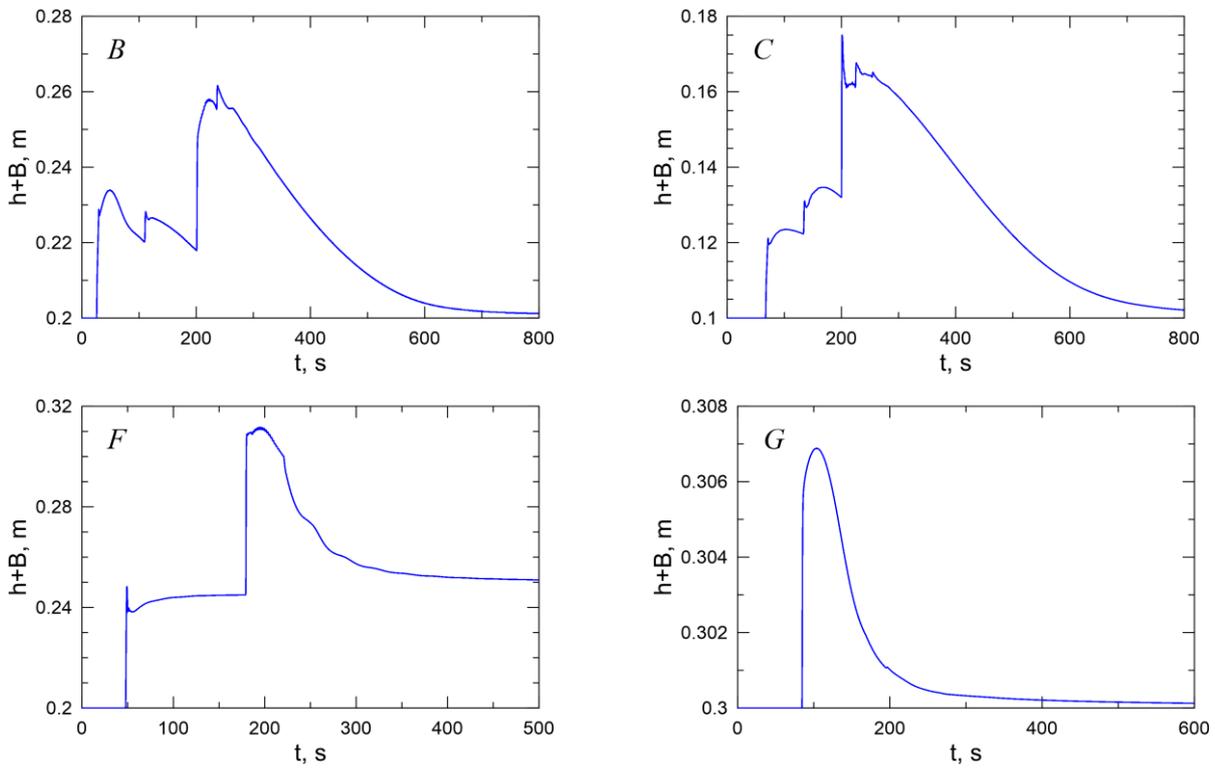

**Fig.34.** Time series of water surface elevations at the junction nodes.

Open water surface elevations along various channels of the network at different times are shown in Fig.29-33. Time serious of water surface elevations at the junction nodes *B*, *C*, *F* and *G* are given in Fig.34. In Fig.29-30 anyone can see that during passing first wave along channel *EF* inflow discharge of water to the junction node *F* does not equal to outflow discharges from it which are zero in this time. Thus, the condition that the sum of the discharges has to be zero at the junction used in many hydrological models as internal boundary condition for junction treatment is not valid in general case.

**Conclusions**

A central-upwind scheme has been applied to simulate shallow water flow in a multiply-connected channel network with arbitrary channel geometry and bottom topography. New reconstruction algorithm of water surface elevation for partially flooded cells has been proposed. This reconstruction algorithm that is a generalization of the idea from [25], with the exact integration of source terms of the shallow water equations provide the well-balanced property and positivity preserving of the scheme. We considered two models of a channel junction treatment based on: 1) continuity equation for a subcritical flow and 2) mass and momentum conservation equations for a supercritical flow. It is shown that for a subcritical flow the continuity equation at a channel junction is generalization of the equality of water surface elevations model. Applying the local draining time approach from [35] to limit outflowing flux from draining cell, we provide the positivity preserving without of a reduction of the CFL time step restriction. Implicit treatment of only a part of the friction slope according to Chertok et al. [39] holds the scheme stability without additional time step restriction.

The set of the numerical tests demonstrates the scheme accuracy, positivity preserving and well-balancing, convergence of numerical results to steady-state solutions and their good agreement with exact solutions and experimental data, including measurements in the multiply-connected river channel network. We propose the new specialized test for simulation of inundation and drying of the channel network by consequence of supercritical and subcritical flows. The test results demonstrate the stability of the scheme and the robustness of the new numerical algorithm for the simulation of the wetting/drying flow in the multiply-connected channel network.